\documentclass[english,12pt]{article}%
\usepackage{bbm,ulem}
\usepackage{graphicx}
\usepackage{subcaption}
\usepackage{makeidx}
\usepackage{amssymb}
\usepackage[utf8]{inputenc}
\usepackage{amsfonts}
\usepackage{amsmath}
\usepackage{amssymb}
\usepackage{amsthm}
\usepackage{url}
\usepackage[colorlinks,citecolor=blue]{hyperref}
\usepackage{mathabx}
\usepackage{wrapfig}
\usepackage{dsfont}
\usepackage{resizegather}
\usepackage{yfonts}
\usepackage{authblk}
\usepackage{geometry}
\usepackage[dvipsnames]{xcolor}

\providecommand{\U}[1]{\protect\rule{.1in}{.1in}}

\usepackage{xcolor}

\newcommand{\hooklongrightarrow}{\lhook\joinrel\longrightarrow}

\newtheorem{prop}{Proposition}[section]
\newtheorem{cor}[prop]{Corollary}
\newtheorem{defi}[prop]{Definition}

\newtheorem{lem}[prop]{Lemma}

\newtheorem{theo}[prop]{Theorem}

\newcommand{\I}{\ensuremath{\mathbb{I}}}
\newcommand{\un}{\mathds{1}}
\newcommand{\tr}{\mbox{\rm Tr}}

\newcommand{\CC}{\mathbb{C}}

\newcommand{\EE}{\mathbb{E}}

\newcommand{\HH}{\mathbb{H}}

\newcommand{\NN}{\mathbb{N}}
\newcommand{\PP}{\mathbb{P}}

\newcommand{\RR}{\mathbb{R}}

\newcommand{\VV}{\mathbb{V}}

\newcommand{\WW}{\mathbb{W}}
\newcommand{\YY}{\mathbb{Y}}
\newcommand{\ZZ}{\mathbb{Z}}
\newcommand{\GG}{\mathbb{G}}

\newcommand{\La}{ {\cal L }}
\newcommand{\Na}{ {\cal N }}

\newcommand{\Ea}{ {\cal E }}
\newcommand{\Sa}{ {\cal S }}

\newcommand{\Va}{ {\cal V }}
\newcommand{\Ua}{ {\cal U }}
\newcommand{\Fa}{ {\cal F }}
\newcommand{\Ga}{ {\cal G }}

\newcommand{\Oa}{ {\cal O }}

\newcommand{\Xa}{ {\cal X }}
\newcommand{\Ma}{ {\cal M }}
\newcommand{\Ta}{ {\cal T}}

\newcommand{\Pa}{ {\cal P }}
\newcommand{\Za}{ {\cal Z }}
\newcommand{\Ya}{ {\cal Y }}
\newcommand{\Wa}{ {\cal W }}

\newcommand{\point}{\mbox{\LARGE .}}
\newcommand{\cqfd}{\hfill\blbx \\}
\def\blbx{\hbox{\vrule height 5pt width 5pt depth 0pt}\medskip}

\def \PP{\mathbb{P}}
\def \RR{\mathbb{R}}

\def \EE{\mathbb{E}}
\def \OO{\mathbb{O}}
\def \CC{\mathbb{C}}

\def \JJ{\mathbb{J}}
\def \ZZ{\mathbb{Z}}
\def \WW{\mathbb{W}}

\newcommand{\vertiii}[1]{{\left\vert\kern-0.25ex\left\vert\kern-0.25ex\left\vert #1
    \right\vert\kern-0.25ex\right\vert\kern-0.25ex\right\vert}}

\begin{document}

% =======================================================
  \title{Fluctuations and Long-Time Stability of Multivariate Ensemble Kalman Filters}

\author{P. Del Moral\thanks{Centre de Recherche Inria Bordeaux Sud-Ouest, Talence, 33405, France. {\footnotesize E-Mail:\,} \texttt{\footnotesize pierre.del-moral@inria.fr}}, B. Nasri \thanks{D\'epartement de m\'edecine sociale et pr\'eventive, \'Ecole de sant\'e publique, Universit\'e de Montr\'eal. {\footnotesize E-Mail:\,} \texttt{\footnotesize bouchra.nasri@umontreal.ca}}
 \& B. R\'emillard \thanks{Department of Statistics and Business Analytics, United Arab Emirates University and Department of Decision Sciences, HEC Montr\'eal. {\footnotesize E-Mail:\,} \texttt{\footnotesize bruno.remillard@hec.ca}}}

\maketitle
  \begin{abstract}
  We develop a self-contained stochastic perturbation theory for discrete-generation and multivariate Ensemble Kalman filters. Unlike their continuous-time counterparts, discrete EnKF algorithms are defined through a two-step prediction–update mechanism and exhibit non-Gaussian fluctuations, even in linear settings. In the multivariate case, these fluctuations take the form of non-central Wishart-type perturbations, which significantly complicate the mathematical analysis.
We establish non-asymptotic, time-uniform stability and error estimates for the ensemble covariance matrix processes under minimal structural assumptions on the signal–observation model, allowing for possibly unstable dynamics. Our results quantify the impact of ensemble size, dimension, and observation noise, and provide explicit bounds on the propagation of stochastic errors over long time horizons.
The analysis relies on a detailed study of stochastic Riccati difference equations driven by matrix-valued non-central Wishart fluctuations. Beyond their relevance to data assimilation, these results contribute to the probabilistic understanding of ensemble-based filtering methods in high dimension and offer new tools for the analysis of interacting particle systems with matrix-valued dynamics.\\



\end{abstract}
%\tableofcontents
\section{Introduction}

\subsection{Ensemble Kalman filters}

The Ensemble Kalman filter (EnKF) is a class of interacting particle system
methods for solving nonlinear filtering and inverse problems, originally
introduced by Evensen in the seminal article~\cite{evensen-intro}; see
also~\cite{evensen-review,evensen-book} for more recent overviews. Over the
past three decades, EnKF methodologies have become among the most widely
used numerical techniques for high-dimensional forecasting and data
assimilation problems, with prominent applications in oceanic and
atmospheric sciences~\cite{allen,lisa,majda,kalnay,ott}, fluid
mechanics~\cite{beyou,memin-1,memin-2}, image inverse
problems~\cite{beyou-2}, weather
forecasting~\cite{anderson-jl,anderson-jl-2,burgers,houte}, environmental
and ecological statistics~\cite{eknes,johns}, oil reservoir
simulation~\cite{evensen-reservoir,nydal,seiler,skj,weng}, and many other
fields.

Connections between EnKF techniques and particle filtering methods for
high-dimensional problems arising in fluid mechanics are discussed
in~\cite{kantas}. We also refer to Section~\ref{sec-EnKF-mf} of the present
article, which is devoted to nonlinear Kalman-type Markov chains and their
mean-field particle interpretations.

The mathematical analysis of continuous-time EnKF models was initiated
in~\cite{DelMoral/Tugaut:2016,dm-k-tu} and further developed in a series of
subsequent works~\cite{BishopDelMoralMatricRicc,dm-b-niclas,bishop-19}; see
also the review articles~\cite{bishop-20,reich-stuart}. Extensions to
nonlinear filtering settings have been considered more recently
in~\cite{lange-1,lange-2,lange-3,lange-4,sahani}. In contrast with the
feedback particle filter framework~\cite{prashant-1,prashant-2,prashant-3,prashant-4},
these studies focus on EnKF particle systems that are, in general, not
consistent for nonlinear models and recover the exact filtering
distribution only in the linear–Gaussian case.

A major advantage of continuous-time formulations is that the ensemble
mean and covariance satisfy coupled nonlinear diffusion equations, which
can be analyzed using classical tools from stochastic analysis. In the
multivariate setting, however, the sample covariance evolves according to
a stochastic Riccati equation driven by matrix-valued martingales, whose
analysis already requires refined techniques in matrix-valued stochastic
calculus~\cite{BishopDelMoralMatricRicc,dm-b-niclas}.

The analysis of discrete-generation EnKF schemes is considerably more
delicate. Unlike their continuous-time counterparts, discrete-time generation EnKF
algorithms are defi\-ned through a two-step prediction–update pro\-ce\-dure
(also referred to as forecasting–analysis steps in the data assimilation
literature) rather than a single coupled diffusion
process~\cite{horton-2}. Moreover, the Gaussian structure underlying the
diffusion limits of continuous-time EnKF models is lost. Even in
one-dimensional settings, discrete-generation EnKF dynamics involve
non-central $\chi^2$-type fluctuations, as illustrated in Theorem~3.1 and
Corollary~3.4 of~\cite{horton-2}. In the multivariate case, this suggests
that the EnKF evolution is governed by analogous two-step transitions
driven by non-central Wishart-type matrix fluctuations.

In the present article, we develop a self-contained and comprehensive
stochastic perturbation analysis of discrete-generation, multivariate
Ensemble Kalman filters. Our results address the fluctuations, stability,
and long-time behavior of EnKF schemes, and include time-uniform,
non-asymptotic mean-square error bounds that remain valid for possibly
unstable signal dynamics.

To the best of our knowledge, these results constitute the first
theoretical analysis of this type for multivariate, discrete-generation
EnKF algorithms. Beyond their relevance for data assimilation, the
stochastic Riccati difference equations studied in this work are of
independent interest, as they provide a prototypical example of a new
class of Markov chains driven by non-central Wishart transitions.

Short proofs are presented in the main text, while more technical
arguments are deferred to the appendices.

\subsection{Some basic notation}
This section presents some basic notation and preliminary results necessary for the statement of our main results.

We denote by $\Ma_{d,d_1}=\RR^{d\times d_1}$ the ring of $(d\times d_1)$-matrices with real entries, for some $d,d_1\geq 1$, when $d=d_1=d$ we write $\Ma_{d}$ instead of $\Ma_{d,d}$. When there is no chance of confusion, we also slightly abuse notation and denote by $0$ and $\I$ the null and identity matrices, respectively, in $\Ma_{d}$ for any dimension $d\geq 1$. We write $A^{\prime}$ to denote the transposition of a matrix $A$.

 When $A\in\Sa_r$ we let $\lambda_1(A)\geq \ldots\geq\lambda_{d}(A)$ denote the ordered eigenvalues of $A$. We equip $\Ma_{r}$ with the spectral norm $\Vert A \Vert=\Vert A \Vert_2=\sqrt{\lambda_{1}(AA^{\prime})}$ or the Frobenius norm $\Vert A \Vert=\Vert A \Vert_{\mathrm{Frob}}=\sqrt{\tr(AA^{\prime})}$.

 We equip the set $\Ma_{d}$ with the spectral norm $\Vert A \Vert_2=\sqrt{\lambda_{max}(AA^{\prime})}$ (a.k.a the $2$-norm) where $\lambda_{max}(\cdot)$ denotes the maximal eigenvalue. The minimal eigenvalue is denoted by $\lambda_{min }(\cdot)$.
We also consider the Frobenius norm
$\Vert A \Vert_{\tiny F}:=\sqrt{\tr(AA^{\prime})}$.

Hereafter, points $x$ in the Euclidean space $\RR^d$ are represented by $d$-dimensional column vectors (or, equivalently, by $d \times 1$ matrices). Note that, for $x\in\RR^{d \times 1}=\RR^d$, the Frobenius and the spectral  norm $\Vert x\Vert_F=\sqrt{x^{\prime}x}=\Vert x\Vert_2$ coincides with the Euclidean norm. 
When there is no possible confusion, sometimes
 we use the notation $\Vert\cdot\Vert$ for any equivalent matrix or vector norm.

We also denote by  $\mbox{\rm Spec}(A)\subset\CC$ the set of eigenvalues of a matrix $A$, and by
$\rho(A)$ the spectral radius of a matrix $A$ defined by
$$
\rho(A):=\max{\left\{\vert \lambda\vert~:~\lambda\in \mbox{\rm Spec}(A)\right\}}.
$$

 We also let $\Sa_d\subset \Ma_{d}$ denote the subset of symmetric matrices, and we let
 $ \Sa_d^0\subset\Sa_d$  { be the closed convex cone} of  {positive} semi-definite matrices,  and
its interior $\Sa_d^+\subset \Sa_d^0$  which resumes to {the open subset} of  {positive definite}  matrices.

 We sometimes use the L\" owner partial ordering notation $S_1\geq S_2$ to mean that a symmetric matrix $(S_1-S_2)$ is positive semi-definite (equivalently, $S_2 - S_1$ is negative semi-definite), and $S_1>S_2$ when $(S_1-S_2)$ is positive definite (equivalently, $S_2 - S_1$ is negative definite). Given $S\in \Sa_d^0-\Sa_d^+$ we denote by $S^{1/2}$ a (non-unique) but symmetric square root of $B$ (given by a Cholesky decomposition). When $S\in\Sa_d^+$ we always choose the principal (unique) symmetric square root. We also denote by $A_{\tiny sym}=(A+A^{\prime})/2$ the symmetric part of a given matrix $A\in\Ma_d$.

 The covariance matrix
of a the probability distribution $\eta$ with second absolute moment on $\RR^d$ is defined by the formula:
$$
{\sf cov}_{\eta}:=\int \left(x- \int z~ \eta(dz)\right)\left(x- \int z~ \eta(dz)\right)^{\prime}\eta(dx)
$$

For a given matrix $U=(U^{j}_{i})_{1\leq i\leq d,1\leq j\leq d_1}\in \Ma_{d,d_1}$, consider the column and row decompositions
$$
U=\left(U^1,\ldots,U^{d_1}\right)=\left(
\begin{array}{c}
U_1\\
\vdots\\
U_{d}
\end{array}
\right)\quad\mbox{\rm with}\quad
U_{i}:=(U_{i}^1,\ldots,U_{i}^{d_1})
\quad\mbox{\rm and}\quad
U^{j}:=\left(
\begin{array}{c}
U^{j}_1\\
\vdots\\
U^{j}_{d}
\end{array}
\right).
$$
We also set
$$
M(U):=\left(M(U)^1,\ldots,M(U)^{d_1}\right)\in  \Ma_{d,d_1}
$$
with
\begin{equation}\label{def-Mmu}
M(U)^{i}=m(U):=\frac{1}{d_1}\sum_{1\leq j\leq d_1} U^{j}
\end{equation}
For symmetric matrices $S\in  \Sa_{d}$ sometimes we write
$
S_{i,j}=S_{j,i}
$ instead of $S^i_j=S^j_i$.

Next we recall a couple of rather well-known estimates in matrix theory. For any $(r\times r)$-square matrices $(P,Q)$ by a direct application of Cauchy-Schwarz inequality we have
\begin{equation}\label{form-pq-ref}
\vert\mbox{\rm tr}(PQ)\vert\leq \Vert P\Vert_F~\Vert Q\Vert_F.
\end{equation}
For any $P,Q\in\Sa^+_d$ we also have 
\begin{equation}\label{f30}
 \mbox{\rm tr}\left(P^2\right)\leq \left(\mbox{\rm tr}\left(P\right)\right)^2\leq d~\mbox{\rm tr}\left(P^2\right)
\quad\mbox{\rm
and}\quad
\lambda_{d}(P)~\mbox{\rm tr}\left(Q\right)\leq \mbox{\rm tr}\left(PQ\right)\leq \lambda_{1}(P)~\mbox{\rm tr}\left(Q\right).
\end{equation}
The above inequality is also valid when $Q$ is positive semi-definite and $P$ is symmetric. We check this claim
using an orthogonal diagonalization of $P$ and recalling that $Q$ remains positive semi-definite (thus with non negative
diagonal entries).
%
%
%
%
%Finally, define the optimal matching distance between the spectrum of matrices $A,B\in\Ma_d$ by
%\begin{equation}\label{optimal-match-d}
%d_{\sf opm}\left(\mathrm{Spec}(A),\mathrm{Spec}(B)\right)=\min_{\mathrm{perm(\cdot)}}\,{\max_{1\leq i\leq d}\vert \lambda_i(A)-\lambda_{\mathrm{perm}(i)}(B)\vert}
% \end{equation}
%where the minimum is taken over the set of $d!$ permutations of $\{1,\ldots,d\}$. Recall also the Krause \cite{KRAUSE199473} and Friedland \cite{friedland1982variation} inequalities,
%\begin{equation}\label{krause-ref}
% d\left(\mathrm{Spec}(A),\mathrm{Spec}(B)\right)\vee \vert\mbox{\rm det}(A)-\mbox{\rm det}(B)\vert^{1/d}~\leq~ c\,\left[
% \Vert A\Vert\vee\Vert B\Vert \right]^{1-1/d}\,\Vert A-B\Vert^{1/d}
%\end{equation}
%for any $A,B\in\Ma_{d}$. For any $A,B\in\Sa_d$ we also have \cite{fiedler2008special} the Hoffman-Wielandt inequality
%\begin{equation}\label{hw-ref}
%\sum_{1\leq i\leq d}\left(\lambda_i(A)-\lambda_i(B)\right)^2~\leq~ \Vert A-B\Vert_F^2
%\end{equation}
For any $P,Q\in\Sa_{r}^+$ we have the Ando-Hemmen inequality
\begin{equation}\label{square-root-key-estimate}
\Vert P^{1/2}- Q^{1/2}\Vert \leq \left[\lambda^{1/2}_{min}(P)+\lambda^{1/2}_{min}(Q)\right]^{-1}~\Vert P- Q\Vert
\end{equation}
for any unitary invariant matrix norm $\Vert . \Vert$ (such as the $2$-norm and the Frobenius norm). See for instance Theorem 6.2 on page 135 in~\cite{higham}, as well as Proposition 3.2  in~\cite{hemmen}. For a more thorough discussion on the geometric properties of positive semidefinite matrices and square roots we refer to~\cite{hiriart}.

Throughout this paper, we denote by $c$ as well as $c_{\alpha}(\beta)$  generic constants
that may depends on some parameters $(\alpha,\beta)$, whose values may vary from line to line but they do not depend on the time parameter.

\section{Description of the models}

\subsection{The Kalman filter}

Consider a Markov chain $(X_n,Y_n)\in (\RR^{d}\times\RR^{d_0})$
 defined for any $n\geq 0$ by the recursive relations
\begin{equation}\label{kbuc}
X_{n+1}=A ~X_{n}+W_n\quad \mbox{\rm and}\quad
Y_n=B  ~X_{n}+V_n
\end{equation}
for some conformal matrices
$(A ,B)$ and some $\RR^{d}$ and
$\RR^{d_0}$-valued  independent centered Gaussian random sequences
$W_n $ and $V_n$ with covariance matrices
$
R >0$, and $
R_0  >0$. The initial random variable $X_0$ is and $ \RR^{d}$-valued Gaussian random variable  independent of the sequence $(W_n,V_n)$ with a mean and
covariance matrix denoted by $\widehat{X}^-_0$ and $P_0>0$.

Let $\Ya_n=\sigma\left(Y_k,~k\leq n\right)$ be the filtration generated by the observation process. The optimal filter is defined by the
conditional distribution $\widehat{\eta}_n$ of the signal state $X_n$ given $\Ya_n$. The distribution $\widehat{\eta}_n$ is a Gaussian  distribution $\Na_d\left(\widehat{X}_n,\widehat{P}_n\right)$ on $\RR^d$ with  conditional mean and covariance
$$
\widehat{X}_n:=\EE(X_n~|~\Ya_n)\quad\mbox{\rm and}\quad
\widehat{P}_n:=\EE\left(\left(X_n-\widehat{X}_n\right)\left(X_n-\widehat{X}_n\right)^{\prime}\right).
$$
The optimal one-step predictor is defined by the conditional distribution $\eta_{n}$ of the signal state $X_n$ given $\Ya_{n-1}$. The distribution $\eta_n$ is also a Gaussian  distribution $\Na_d\left(\widehat{X}^-_n,P_n\right)$ on $\RR^d$ with conditional  mean and covariance
$$
\widehat{X}^-_n:=\EE(X_n~|~\Ya_{n-1})\quad\mbox{\rm and}\quad
P_n:=\EE\left(\left(X_{n}-\widehat{X}_{n}^{-}\right)\left(X_{n}-\widehat{X}_{n}^{-}\right)^{\prime}\right).
$$
\subsubsection*{Updating-prediction transitions}
The updating-prediction steps  are described by the following synthetic diagram
$$
(\widehat{X}_n^-,P_n)\longrightarrow (\widehat{X}_{n},\widehat{P}_{n})\longrightarrow (\widehat{X}^-_{n+1},P_{n+1}),
$$
with the well-known updating-prediction Kalman filter equations:
\begin{equation}\label{kalman-rec}
\left\{\begin{array}{rcl}
\widehat{X}_{n}&=&\widehat{X}^{-}_{n}+K(P_n)
~(Y_n-B  \widehat{X}^{-}_{n})  \\
&&\\
 \widehat{P}_{n}&=&\widehat{A}(P_n)~P_{n}
\end{array}\right.\quad\mbox{\rm and} \quad\left\{\begin{array}{rclcrcl}
\widehat{X}^{-}_{n+1}&=&A \widehat{X}_{n}\\
&&\\
P_{n+1}&=&A \widehat{P}_{n}A^{\prime} +R,
\end{array}\right.
\end{equation}
In the above display, $P\mapsto K(P)$ stands for the so-called  Kalman gain function defined for any $P\geq 0$  by the formula
\begin{equation}\label{def-Aw}
 \begin{array}{l}
\displaystyle K(P):=PB^{\prime}(B  PB  ^{\prime}+R_0 )^{-1}
\quad \mbox{\rm and }\quad 
  \widehat{A}(P):=\I-K(P) B
 \end{array}
 \end{equation}

 \subsubsection*{Riccati difference equation}
We have
$$
P-PSP\leq \widehat{A}(P)P=P-PB^{\prime}(B  PB  ^{\prime}+R_0 )^{-1}BP\leq P
$$
 as well as
\begin{equation}\label{eq:Gn}
 \widehat{A}(P)=(\I+PS)^{-1} \quad\mbox{\rm with} \quad S:=B  ^{\prime}R_0^{-1} B.
 \end{equation}
 This yields the off-line Riccati covariance evolution
\begin{equation}\label{def-Riccati-drift}
 P_{n+1}
=\Phi(P_n) 
\end{equation}
with the Riccati map
$$
P\mapsto \Phi(P):=A ~(\I+P S)^{-1}~ P~A ^{\prime}+R
$$
Moreover, we have
$$
R\leq
\Phi(P)=APA^{\prime}-APB^{\prime}(R_0+BPB^{\prime})^{-1}BPA^{\prime}+R
\leq APA^{\prime}+R
$$
as well as
$$
\Phi(P):=A ~(\I+P S)^{-1}~ P~A ^{\prime}+R\leq A ~S^{-1}~A ^{\prime}+R
$$

We also quote the symmetric formulation of the Riccati updating map
\begin{equation}\label{eq:gnB-ident1}
\widehat{A}(P)P=\widehat{A}(P)P\widehat{A}(P)^{\prime}+\widehat{R}(P)
\quad \mbox{\rm with}\quad  \widehat{R}(P):=K(P)\,R_0\,K(P)^{\prime}
\end{equation}
Moreover we have  the square root formula
\begin{eqnarray}
\widehat{A}(P)P&=&P^{1/2}\left(\I-P^{1/2}B^{\prime}(B  PB  ^{\prime}+R_0 )^{-1}BP^{1/2}\right)P^{1/2}\nonumber\\
&=&P^{1/2}~\left(\I+P^{1/2}SP^{1/2}\right)^{-1}P^{1/2}\label{square-root}
\end{eqnarray}

The proofs of (\ref{eq:Gn}), (\ref{eq:gnB-ident1}), (\ref{square-root}) and (\ref{def-Riccati-drift})  follow elementary matrix product manipulations, thus there are provided in the appendix on page~\pageref{def-Riccati-drift-proof}.

\subsection{The Ensemble Kalman filter}\label{sec-EnKF-mf}
\subsubsection*{Nonlinear Markov chain}
Let $(\Xa_0, \Wa_n,\Va_n)$ be independent copies of $(X_0, W_n,V_n)$.  Consider the nonlinear Markov chain starting at $\Xa_0$ and defined sequentially for any $n\geq 0$ by the updating-prediction formulae
\begin{equation}\label{kalman-X-goth}
\left\{
\begin{array}{rcl}
\widehat{\Xa}_{n}&=&\Xa_n+K({\sf cov}_{\eta_n})~(Y_n-(B\,\Xa_n+\,\Va_n))\quad \mbox{\rm with}\quad \eta_n:=\mbox{\rm Law}(\Xa_n~|~\Ya_{n-1})\\
&&\\
\Xa_{n+1}&=& A\,\widehat{\Xa}_{n}+\Wa_{n+1}.\\
\end{array}
\right.
\end{equation}
Observe that the above stochastic model can be interpreted as a Markov chain
with a two-step updating-prediction transition  described by the following synthetic diagram
$$
\Xa_n\longrightarrow \widehat{\Xa}_{n}\longrightarrow \Xa_{n+1}.$$
The updating transition $\Xa_n\longrightarrow \widehat{\Xa}_{n}$
depends on the observation $Y_n$ delivered by the sensor as well as on
the conditional law of $\Xa_n$ given the information $\Ya_{n-1}$.

Using a simple induction argument, it is straightforward to show that
$$
\begin{array}{rcl}
\eta_n&=&\Na_d\left(\widehat{X}^-_n,P_n\right)=\mbox{\rm Law}(\Xa_n~|~\Ya_{n-1})\quad\Longrightarrow\quad {\sf cov}_{\eta_n}=P_n,\\
\\
 \widehat{\eta}_n&=&\Na_d\left(\widehat{X}_n,\widehat{P}_n\right)=\quad\mbox{\rm Law}(\widehat{\Xa}_n~|~\Ya_{n}) \quad\Longrightarrow\quad {\sf cov}_{ \widehat{\eta}_n}= \widehat{P}_n.
\end{array}
$$

Rewritten in a slightly different form the updating transition in (\ref{kalman-X-goth}) takes the form
\begin{equation}\label{def-w-AR}
\begin{array}{l}
\widehat{\Xa}_{n}=\widehat{A}(P_n)~\Xa_n+K(P_n)~Y_n+\widehat{\Va}_n\quad
 \mbox{\rm with}\quad
\widehat{\Va}_n\sim\Na_d(0,\widehat{R}(P_n))
\end{array}
\end{equation}

\subsubsection*{Mean field particle algorithm}
The Ensemble Kalman filter (abbreviated EnKF) associated to the filtering problem (\ref{kbuc}) coincides with the mean-field particle interpretation of the nonlinear Markov chain (\ref{kalman-X-goth}).
To define these models, we let $\xi_0=(\xi_0^i)_{1\leq i\leq N+1}$ be a sequence of $(N+1)$ independent copies of  $X_0$, for some parameter $N\geq 1$. We also denote by $(\Wa^i_n)_{i\geq 1}$ and $(\Va^i_n)_{i\geq 1}$ a sequence of independent copies of $W_n$ and $V_n$.

 The mean field particle interpretation of the nonlinear Markov chain discussed above is given by the interacting particle system defined sequentially for any $1\leq i\leq N+1$ and $n\geq 0$ by the formulae
 \begin{equation}\label{kalman-EnKF-def}
\widehat{\xi}^i_{n}=\xi^i_n+K(p_n)~(Y_n-(B\xi^i_n+\Va^i_n))
\quad \mbox{\rm and}\quad
\xi^i_{n+1}=A\,\widehat{\xi}^i_{n}+\Wa^i_{n+1}.
\end{equation}
with the normalised sample covariance $p_n$ given by
\begin{equation}\label{scov-kalman-EnKF}
p_n:=\frac{1}{N}\sum_{1\leq i\leq N+1}(\xi^i_n-m_n)(\xi^i_n-m_n)^{\prime}=\left(1+\frac{1}{N}\right) {\sf cov}_{\eta^N_n}.
\end{equation}
In the above display, $\eta^N_n$ and $ m_n$ stand for the empirical measure and the sample mean defined by
$$
\eta^N_n:=\frac{1}{N+1}\sum_{1\leq i\leq N+1}\delta_{\xi^i_n}\quad \mbox{\rm and}\quad
 m_n:=\frac{1}{N+1}\sum_{1\leq i\leq N+1}\xi^i_n.
$$
We also consider the updated sample means and sample covariance matrices
\begin{equation}\label{sm-kalman-EnKF}
 \widehat{m}_n:=\frac{1}{N+1}\sum_{1\leq i\leq N+1}\widehat{\xi}^i_n
 \quad\mbox{\rm and}\quad \widehat{p}_n:=\frac{1}{N}\sum_{1\leq i\leq N+1}(\widehat{\xi}^i_n-\widehat{m}_n)(\widehat{\xi}^i_n-\widehat{m}_n)^{\prime}
\end{equation}

Whenever $B=0$ is the null matrix we have $K(p_n)=0$ so that $\widehat{\xi}^i_n=\xi^i_n$ and therefore  $\widehat{p}_n=p_n$. In this situation the sample covariance $p_n$ resume to the one associated with independent copies of a linear Gaussian Markov chain
$$
\xi^i_{n+1}=A\xi^i_n+\Wa^i_{n+1}
$$
In this particular case, the random matrices $p_n$ are distributed according to a Wishart distribution.

We set
$$
\xi_n:=\left(\xi_n^1,\ldots,\xi_n^{N+1}\right)\quad \mbox{\rm and}\quad
\widehat{\xi}_n:=\left(\widehat{\xi}_n^{\,1},\ldots,\widehat{\xi}_n^{\,N+1}\right)
$$
We also consider the matrices
$$
\WW_n:=\left[\Wa_n^1,\ldots,\Wa_n^{N+1}\right]\in \Ma_{d\times(N+1)}
\quad \mbox{\rm and}\quad
\VV_n:=\left[\Va_n^1,\ldots,\Va_n^{N+1}\right]\in \Ma_{d_0\times(N+1)}
$$
and we set
$$
\YY_n:=\left(Y_n,\ldots,Y_n\right)\in \Ma_{d_0\times(N+1)} 
$$
In this notation, the updating and the prediction of the EnKF take the following matrix form
\begin{eqnarray}
\widehat{\xi}_{n}&=&\xi_n+K(p_n)~(\YY_n-(B\xi_n+\VV_n))
\nonumber\\
&&\nonumber\\
\xi_{n+1}&=&A ~\widehat{\xi}_{n}+ \WW_{n+1}\label{matrix-form}\end{eqnarray}

\subsection{Some time-uniform estimates}

The next theorem provides time-uniform mean-error  for  the sample covariance matrices $p_n$ defined in (\ref{scov-kalman-EnKF}).
\begin{theo}\label{th-intro-punif}
For any $n\geq 1$ we have the under-bias property
\begin{equation}\label{unif-ulb}
0<R\leq \EE(p_{n})\leq P_n\leq  AS^{-1}A^{\prime}+R
\end{equation}
There exists some parameter $\iota>0$ such that for any $n\geq 0$ and $N>(1+d)$ we have
$$
0\leq P_n-\EE(p_n)\leq \frac{\iota}{N}~\I
$$
In addition,
 for any $r\geq 1$ and $N>(1+d)$ and we have the time-uniform estimates
\begin{equation}\label{unif-cov}
\sup_{n\geq 0}\EE\left(\Vert\, p_n-P_n\Vert^r\right)^{1/r}\leq c(r)/\sqrt{N}
\end{equation}
\end{theo}
The detailed proof of the above theorem is provided on page~\pageref{th-intro-punif-proof}.

The next corollary provides time-uniform mean-error  estimates for the updated sample covariance and the sample gain matrix of the EnKF filter.
\begin{cor}\label{cor-wpg}
For any $r\geq 1$ and $N+1>d$ and we have the time-uniform estimates
\begin{equation}\label{unif-cov-2}
\sup_{n\geq 0}\left(\EE\left(\Vert\, \widehat{p}_n-\widehat{P}_n\Vert^r\right)^{1/r}\vee \EE\left(\Vert\, K(p_n)-K(P_n)\Vert^r\right)^{1/r}\right)\leq c(r)/\sqrt{N}
\end{equation}
\end{cor}
The detailed proof of the above corollary is provided in the appendix, on page~\pageref{cor-wpg-proof}.

\section{Riccati difference equations}

This section is devoted to the analysis of discrete-time Riccati difference
equations and to the matrix-valued objects that naturally arise in their
perturbation theory. These equations play a central role in the stability
and fluctuation analysis of Ensemble Kalman filters, as they govern the
evolution of the empirical covariance matrices and their first order linearizations.

Our main objective is to collect and develop a set of analytical tools
that will be repeatedly used throughout the paper. In particular, we
establish differentiability properties, contraction estimates, and
explicit decomposition formulae for Riccati semigroups and their associated
matrix products. Although many of these results are classical in control
theory, we present them here in a unified and self-contained form adapted
to the probabilistic setting considered in this work.

The section is organized as follows. In Section~\ref{sec-diff-calculus},
we recall basic elements of Fr\'echet differential calculus for
matrix-valued maps on the cone of symmetric matrices, which are required
to control first- and second-order perturbations of Riccati operators.
Section~\ref{sec-reg-props} reviews fundamental regularity and monotonicity
properties of discrete Riccati semigroups, together with their stability
behavior around the unique positive definite fixed point. Finally,
Section~\ref{sec-floquet} presents a Floquet-type representation for
directed matrix products associated with Riccati flows, leading to
explicit first- and second-order expansion formulae that are crucial for
the non-asymptotic error analysis developed in subsequent sections.

\subsection{Differential calculus}\label{sec-diff-calculus}
We let $\La(\Sa_d,\Ma_{d})$ be the set of bounded linear  functions from $\Sa_d$ into $\Ma_d$. Let  {$\Oa_d\subset\Ma_d$ be a non empty open and convex subset of $\Sa_d$}.
We recall that a mapping  {$\Upsilon:\Oa_{d}\mapsto \Sa_d$} defined in some domain $\Oa_{d}$
is Fr\'echet differentiable at some $A\in \Oa_d$ if there exists a
continuous linear  function
$
\nabla \Upsilon(A)\in  \La(\Sa_d,\Ma_{d})
$
such that
$$
\lim_{\Vert H\Vert\rightarrow 0}\Vert H\Vert^{-1}
\Vert \Upsilon(A+H)-\Upsilon(A)-\nabla \Upsilon(A)\cdot H \Vert=0
$$
In other words, for any given $A\in \Oa_d$ and $\epsilon>0$ there exists some $\delta>0$ such that
$$
\Vert H\Vert\leq \delta\Longrightarrow A+H\in \Oa_d\quad\mbox{\rm and}\quad
\Vert \Upsilon(A+H)-\Upsilon(A)-\nabla \Upsilon(A)\cdot H \Vert\leq \epsilon~\Vert H\Vert
$$
The l.h.s. condition is met for $\Oa_d=\Sa^+_d\ni A$. We check this claim using Weyl's inequality
$$
\lambda_{ min}(A+H)\geq \lambda_{ min}(A)+\lambda_{ min}(H)\geq  \lambda_{ min}(A)-\Vert H\Vert_2
$$
This shows that
$$
\Vert H\Vert_2< \lambda_{ min}(A)\Longrightarrow A+H\in \Sa^+_d
$$
 {The function $\Upsilon$ is said to be Fr\'echet  differentiable on $\Oa_d$ when the mapping $$\nabla \Upsilon~:~A\in \Oa_d
\mapsto \nabla \Upsilon(A)\in  \La(\Sa_d,\Ma_{d})$$ is continuous.
Higher Fr\'echet derivatives are defined in a similar way. For instance, the mapping $\Upsilon$  is twice Fr\'echet differentiable at $A\in \Oa_d$ when
the mapping
$
\nabla \Upsilon
$
is also Fr\'echet  differentiable at $A\in \Oa_d$}. Identifying $ \La(\Sa_d,\La(\Sa_d,\Ma_{d}))$ with the set
$ \La(\Sa_d\times\Sa_d,\Ma_{d})$ of continuous bilinear maps from $(\Sa_d\times\Sa_d)$ into $\Sa_d$, the second derivative
$$\nabla^2 \Upsilon~:~A \in \Oa_d \mapsto \nabla^2 \Upsilon(A)\in  \La(\Sa_d\times\Sa_d,\Ma_{d})$$ is defined by a continuous and symmetric billinear map $\nabla^2 \Upsilon(A)$ such that the limit
$$
\lim_{\Vert H_2\Vert\rightarrow 0}\Vert H_2\Vert^{-1}
\Vert \nabla \Upsilon(A+H_2)\cdot H_1-\nabla \Upsilon(A)\cdot H_1-\nabla^2 \Upsilon(A)\cdot (H_1,H_2) \Vert=0
$$
exists uniformly w.r.t. $H_1\in\Sa_d$ in bounded sets. The polarization formula
$$
\begin{array}{l}
\nabla^2 \Upsilon(A)\cdot (H_1,H_2)\\
\\
\displaystyle=\frac{1}{4}~\left[
\nabla^2 \Upsilon(A)\cdot (H_1+H_2,H_1+H_2)-\nabla^2 \Upsilon(A)\cdot (H_1-H_2,H_1-H_2)\right]
\end{array}
$$
shows that it suffices to compute the second order derivatives $\nabla^2 \Upsilon(A)\cdot (H,H)$ in the same direction $H=H_1=H_2$. Identifying
$\La(\Sa_d\times\Sa_d,\Ma_{d})$ with $\La(\Sa_d\otimes\Sa_d,\Ma_{d})$ sometimes we set $\nabla^2 \Upsilon(A)\cdot H^{\otimes 2}$ instead of $\nabla^2 \Upsilon(A)\cdot (H,H)$. For a more detailed discussion on these tensor product identifications of symmetric multilinear maps we refer to chapter 5 in~\cite{dudley}, see also~\cite{niclas}.
\subsection{Some regularity properties}\label{sec-reg-props}

In this section, we discuss some of the theory behind discrete time Riccati matrix evolution equations and present some results that are of use throughout the paper. This section is mainly taken from~\cite{horton}.

Consider the matrix map $\Ea$ and $\Fa$ defined for any $P\in \Sa^0_{d}$ by the formulae
\begin{equation}\label{E-F-def}
\Ea(P):=A(\I+PS)^{-1}\in \Ma_{d}
\quad \mbox{\rm and}\quad
\Fa(P):=S(\I+PS)^{-1}\in \Sa^0_{d}
\end{equation}
By Lemma 2.1 in~\cite{horton},
for any $P\in\Sa_d^0$ we have
\begin{equation}\label{E-F-estim}
\alpha_-(P)S\leq \Fa(P)=S^{1/2}(\I+S^{1/2}PS^{1/2})^{-1}S^{1/2}\leq \alpha_+(P)S,
\end{equation}
with the positive parameters
$$
\alpha_-(P):=(1+\lambda_{1}(P)\lambda_{1}(S))^{-1}~\quad \mbox{and}\quad
\alpha_+(P):=(1+\lambda_{d}(P)\lambda_{d}(S))^{-1}.
$$

We denote by $\Phi_n$  the Riccati evolution semigroups defined by the inductive composition formula
$\Phi_n=\Phi\circ\Phi_{n-1}$ with the convention $\Phi_0=\I$, the identity map from $\Sa^0_d$ into itself.
Whenever $R>0$ and $R_0>0$ the Riccati evolution semigroups $\Phi_n$   is positive definite, in the sense that
\begin{equation}\label{ref-positive}
P\geq 0\Longrightarrow \forall n\geq 1\qquad \Phi_n(P)>0.
\end{equation}
In addition,  $\Phi$ has unique positive definite fixed point  $\Phi(P_{\infty})=P_{\infty}\in \Sa^+_{d}$
\begin{equation}\label{def-B}
\Ea(P_{\infty}):=A(\I+P_{\infty}S)^{-1}
\quad\mbox{\rm satisfies} \quad \rho(\Ea(P_{\infty}))<1.
\end{equation}
The proof of these assertions under weaker observability-controllability conditions
can be found in any textbook on Riccati equations, see for instance~\cite{kailath,kucera,lancaster} the more recent book~\cite{hisham} and the article~\cite{horton}.  In optimal control theory, the matrix $\Ea(P_{\infty})$ is often called the closed loop-matrix.

 For the convenience of the reader, we recall some rather well-known properties of the Riccati semigroup.
\begin{lem}[\cite{horton}]\label{lem-phi-props}
For any $P,Q\in \Sa^0_{d}$ we have the formulae
\begin{equation}\label{form-1}
\Ea(Q)=\Ea(P)~\left(\I+(P-Q)\Fa(Q)\right)
\quad
\mbox{and}
\quad
\Phi(P)-\Phi(Q)=\Ea(P)(P-Q)\Ea(Q)^{\prime}.
\end{equation}
In addition, the Riccati map $\Phi$ introduced in~\eqref{def-Riccati-drift} satisfies the following properties.
\begin{enumerate}
\item[(i)] For all $P \in \Sa_d^0$, $\Phi(P)' = \Phi(P)$.
\item[(ii)] For $P, Q \in \Sa_d^0$ and $n \ge 1$, $P \ge Q \Longrightarrow \Phi_n(P) \ge \Phi_n(Q)$ and $\Phi_{n+1}(0) \ge \Phi_n(0)$.
\end{enumerate}
\end{lem}

We now consider the directed matrix product $\Ea_n(P)$ defined by
\begin{equation}\label{En-def}
\Ea_{n+1}(P):=\Ea_{n}(\Phi(P))\,\Ea(P)\quad\mbox{\rm with}\quad
\Ea_0(P)=\I\quad\Longrightarrow \quad \Ea_n(P_{\infty})=\Ea(P_{\infty})^n.
\end{equation}

Note that this implies that $\Ea_1 = \Ea$. The matrices $\Ea_n(P)$ play a crucial role in the regularity analysis and the stability theory of Riccati difference equations.
For instance, for any $n\geq 0$ and any $P,Q\in \Sa^0_{d}$ we have the well-known formula
\begin{equation}\label{phi-form}
\Phi_n(P)-\Phi_n(Q)=\Ea_n(P)\,(P-Q)\,\Ea_n(Q)^{\prime},
\end{equation}
whose proof is a simple consequence of \eqref{form-1}. Observe that
 $\Phi_n$ is a smooth matrix functional with a first order Fr\'echet derivative (see~\cite{niclas}) defined for any $P\in\Sa_{d}^0$ and $H\in\Sa_{d}$ by
\begin{equation}\label{fund-eq}
\begin{array}{l}
\displaystyle\nabla\Phi_n(P)\cdot H\,=\,\Ea_n(P)\,H\,\Ea_n(P)^{\prime}\\
\\
\Longleftrightarrow\quad \forall n\geq 1\quad
\nabla\Phi_n(P)=\nabla\Phi_{n-1}(\phi(P))\circ \nabla\Phi(P)\quad \mbox{\rm with}\quad
\nabla\Phi_0(P):=\I.
\end{array}
\end{equation}
In the above display, the symbol $``\circ"$ stands for the composition of operators.

From this perspective, the directed matrix product $\Ea_n(P)$ can be seen as  a fundamental solution of
the first variational equation \eqref{fund-eq} associated with the Riccati difference equation.
Moreover, setting
$$
P_n=\Phi(P_{n-1})\quad \mbox{\rm and}\quad Q=P_{\infty},
$$
the formula \eqref{phi-form} yields the product formula
\begin{equation}\label{phi-form-infty}
\displaystyle P_n-P_{\infty}=\Ea_n(P_0)\,(P_0-P_{\infty})\,(\Ea(P_{\infty})^n)^{\prime},
\end{equation}
where we may write $$\Ea_{n}(P_0):=\Ea(P_{n-1})\ldots  \Ea(P_{1})\, \Ea(P_{0}).
$$

 The spectral radius $\rho(\Ea(P_{\infty}))$ is connected to any norm $\Vert\point\Vert$ of the matrix powers $\Ea(P_{\infty})^n$ arising in \eqref{phi-form-infty} by
 Gelfand's formula (see for instance~\cite{rota}) given by
\begin{equation}\label{gelfand-form}
\rho(\Ea(P_{\infty}))=\lim_{k\rightarrow \infty}\Vert \Ea(P_{\infty})^k\Vert^{1/k}<1,
\end{equation}
where the latter inequality holds due to \eqref{def-B}. Thus, our observability and controllability conditions ensure the exponential decays of  the matrix norms $\Vert \Ea(P_{\infty})^n\Vert$ towards $0$ for sufficiently large time horizons.
\subsection{A Floquet-type formula}\label{sec-floquet}
Consider the Grammian
$$
\Ga_n=\sum_{0\leq k<n}\left(\Ea(P_{\infty})^k\right)^{\prime}\Fa(P_{\infty})~\Ea(P_{\infty})^k\in \Sa^0_d
$$
and let $\La_n$ be the matrix map defined for any $P\in  \Sa^0_d$ by
$$
\begin{array}{l}
\displaystyle
\La_n(P):=
\I+(P-P_{\infty})~\Ga_n\in\Ma_d\\
\\
\displaystyle\Longrightarrow\quad \La_n(P_{\infty})=\I
\quad \mbox{and}\quad
\La_n(P)-\La_n(Q)=(P-Q)\,\Ga_n.
\end{array}$$
Using \eqref{E-F-estim}, for any $n\geq 1$ we check that
$$
\Ga_n\geq \alpha_-(P_{\infty})~\sum_{0\leq k<n}\left(\Ea(P_{\infty})^k\right)^{\prime}S~\Ea(P_{\infty})^k\geq \alpha_-(P_{\infty})~S
$$
In addition, applying the duality formula presented in Theorem 1.1 in~\cite{horton}, we have
$$
\Ga_n\leq\Ga:=\sum_{k\geq 0}\left(\Ea(P_{\infty})^k\right)^{\prime}\Fa(P_{\infty})~\Ea(P_{\infty})^k=\left(P_{\infty}+\overline{P}_{\infty}^{-1}\right)^{-1}< P_{\infty}^{-1}
$$
In the above display $\overline{P}_{\infty}$ stands for the positive definite fixed point of the mapping $\overline{\Phi}$ defined as $\Phi$ by replacing $(A,R,S)$ by $(A^{\prime},S,R)$.

Observe that
\begin{equation}\label{ref-Lan}
\La_n(P)^{-1}=\Ga^{-1}_n~\left(P+
\left(\Ga_n^{-1}-P_{\infty}\right)\right)^{-1}~\quad \mbox{\rm and}\quad
\Ga_n^{-1}-P_{\infty}>0
\end{equation}
with
$$
\Ga^{-1}_n\leq \frac{1}{\alpha_-(P_{\infty})}~S^{-1}\quad \mbox{\rm and}\quad
\left(P+
\left(\Ga_n^{-1}-P_{\infty}\right)\right)^{-1}\leq \left(P+\overline{P}_{\infty}^{-1}\right)^{-1}\leq \overline{P}_{\infty}
$$
Note that
$$
\Ga_1=\Fa(P_{\infty})=
(S^{-1}+P_{\infty})^{-1}\quad\mbox{\rm and}\quad
\La_1(P)^{-1}=(S^{-1}+P_{\infty})~(P+S^{-1})^{-1}
$$
and therefore
$$
\Ea(P_{\infty})~\La_1(P)^{-1}=AS^{-1}(P+S^{-1})^{-1}=\Ea(P)
$$
More interestingly, by theorem  1.3 in~\cite{horton} (see also formula (16) in~\cite{orfanidis})
for any $P\in \Sa_d^0$ we have the Floquet-type Riccati matrix product formula
\begin{equation}\label{form-Floquet}
\Ea_n(P)=\Ea(P_{\infty})^n~\La_n(P)^{-1}\quad \mbox{with}\quad \iota:=
\sup_{P\in \Sa^0_r}\sup_{n\geq 1}\Vert \La_n(P)^{-1}\Vert<\infty.
\end{equation}
This yields
$$
\begin{array}{l}
\left(\La_n(Q)^{-1}-\La_n(P)^{-1}\right)\\
\\
=\La_n(Q)^{-1}\left(\La_n(P)-\La_n(Q)\right)
\La_n(P)^{-1}\\
\\
=\La_n(P)^{-1}~\left(P-Q\right)~\Ga_n~
\La_n(P)^{-1}+
\left(\La_n(Q)^{-1}-\La_n(P)^{-1}\right)~\left(P-Q\right)~\Ga_n~
\La_n(P)^{-1}
\end{array}$$
from which we check the second order decomposition
$$
\begin{array}{l}
\left(\La_n(Q)^{-1}-\La_n(P)^{-1}\right)\\
\\
=\La_n(P)^{-1}\left(P-Q\right)~\Ga_n~
\La_n(P)^{-1}\\
\\
\hskip3cm+
\La_n(Q)^{-1}~\left(P-Q\right)~\Ga_n~
\La_n(P)^{-1}~\left(P-Q\right)~\Ga_n~
\La_n(P)^{-1}
\end{array}$$
Rewritten in a slightly different form we have
\begin{eqnarray*}
\left(\La_n(Q)^{-1}-\La_n(P)^{-1}\right)&=&\nabla \La_n(Q,P)^{-1}\cdot (Q-P)\\
&=&\nabla \La_n(P)^{-1}\cdot (Q-P)+\nabla^2 \La_n(Q,P)^{-1}\cdot (Q-P,Q-P)
\end{eqnarray*}
with
$$
\nabla \La_n(Q,P)^{-1}\cdot H:=-~\La_n(Q)^{-1}~H~\Ga_n~
\La_n(P)^{-1}\quad\mbox{\rm and}\quad
\nabla \La_n(P)^{-1}:=\nabla \La_n(P,P)^{-1}$$
and the second order remainder term
$$
\nabla^2 \La_n(Q,P)^{-1}\cdot (H,H):=
\La_n(Q)^{-1}~H~\Ga_n~
\La_n(P)^{-1}~H~\Ga_n~
\La_n(P)^{-1}
$$
This yields
$$
\begin{array}{l}
\Ea_n(Q)-\Ea_n(P)\\
\\
=\Ea_n(P)~\left(P-Q\right)~\Ga_n~
\La_n(P)^{-1}\\
\\
\hskip3cm+\Ea_n(Q)~\left(P-Q\right)~~\Ga_n~
\La_n(P)^{-1}~\left(P-Q\right)~\Ga_n~
\La_n(P)^{-1}
\end{array}
$$
Rewritten in a slightly different form we have
\begin{prop}\label{prop-dec-Ea}
For any $P,Q\in\Sa^0_d$ and $n\geq 1$ we have the first order decomposition
$$
\Ea_n(Q)-\Ea_n(P)
=\nabla \Ea_n(Q,P)\cdot (Q-P)
$$
with
$$
\nabla \Ea_n(Q,P)\cdot H:=-\Ea_n(Q)~H~\Ga_n~
\La_n(P)^{-1}
$$
We also have the second order decomposition
$$
\Ea_n(Q)-\Ea_n(P)
=\nabla \Ea_n(P)\cdot (Q-P)+\nabla^2 \Ea_n(Q,P)\cdot (Q-P,Q-P)
$$
with
$
\nabla \Ea_n(P):=\nabla \Ea_n(P,P)
$
and the second order remainder
$$
\nabla^2 \Ea_n(Q,P)\cdot (H,H):=\Ea_n(Q)~H~\Ga_n~
\La_n(P)^{-1}~H~\Ga_n~
\La_n(P)^{-1}
$$
\end{prop}

\section{Non central Wishart matrices}

This section is devoted to the study of non-central Wishart random matrices
and their fluctuation properties. Such matrices arise naturally in the
analysis of discrete-generation Ensemble Kalman filters, where sample
covariance matrices are formed from Gaussian ensembles with nonzero means
induced by prediction and update steps. In contrast with the continuous-time
setting, these covariance fluctuations are intrinsically non-Gaussian and
require a careful matrix-valued probabilistic analysis.

Our primary goal is to derive explicit decompositions and moment estimates
for non-central Wishart matrices that are suitable for non-asymptotic and
time-uniform error analysis. Rather than relying on exact distributional
formulae—which are typically expressed in terms of zonal polynomials or
hypergeometric functions of matrix arguments—we adopt an approach based on
orthogonal decompositions, perturbative expansions, and concentration
inequalities. This framework provides transparent bounds that scale
explicitly with the dimension and the ensemble size.

The section is organized as follows. In Section~\ref{sec-gaussian-matrices},
we recall basic properties of Gaussian random matrices and introduce the
non-central Wishart distributions that appear throughout the paper.
Section~\ref{sec-orth-decomp} establishes orthogonal decompositions for
sample covariance matrices, separating deterministic bias terms from
stochastic fluctuations. Section~\ref{sec-fluct-analysis} is devoted to the
analysis of these fluctuations, including multivariate central limit
theorems for the associated matrix-valued noise. Finally,
Section~\ref{sec-inv-moments} provides uniform bounds on inverse moments of
non-central Wishart matrices, which play a crucial role in controlling the
stability of discrete Riccati updates and Ensemble Kalman filter covariance
flows.

\subsection{Gaussian matrices}\label{sec-gaussian-matrices}

We use notation $Z\sim \eta$ for saying that a random variable $Z$ has distribution $\eta$.

Let
$Z=(Z^1,\ldots,Z^N)\in \Ma_{d\times N}$ be a random matrix with $N$
independent and centered Gaussian column $$Z^i=\left(
\begin{array}{c}
Z^{i}_1\\
\vdots\\
Z^{i}_{d}
\end{array}
\right)
\sim \Na_d(0,R)\quad\mbox{\rm with} \quad R\in \Sa_d^0.$$
Consider the covariance mapping
$$
 q~:~z\in \Ma_{d\times N}\mapsto q(z)=\frac{1}{N}~zz^{\prime}\in \Sa^0_{d}
$$
and a matrix square root function
$$
 q^{1/2}~:~z\in \Ma_{d\times N}\mapsto q^{1/2}(z):=q(z)^{1/2}\in\Sa^0_{d} \quad \mbox{\rm with entries}\quad
q^{1/2}_{i,j}(z):=\left(q(z)^{1/2}\right)_{i,j}
$$
\begin{defi}
The  non-central Wishart distribution  $ \mbox{\it Wishart}_d\left(N, Nq(z),R\right)$ on $\Sa^0_d$ with parameters $(N,N q(z),R)$     is the distribution of the random  matrix defined below
$$
\begin{array}{l}
\displaystyle N~q(z+Z)=\sum_{1\leq i\leq N}(z^i+Z^i)(z^i+Z^i)^{\prime}
\end{array}$$
The parameters $R$ and $Nq(z)$ are called respectively the scale matrix (a.k.a. parameter matrix) and the non-centrality parameters. \end{defi}
Sometimes  the non centrality parameter is represented by the matrix $$\frac{N}{2}~R^{-1/2}q(z)R^{-1/2}\quad\mbox{\rm or by the matrix}\quad NR^{-1}q(z),$$
see for instance~\cite{kourouklis} and respectively~\cite{gupta}. 
Central Wishart matrices $Nq(Z)$ correspond to the case $z=0$. In this situation, we write $ \mbox{\it Wishart}_d\left(N,R\right)$ instead of $ \mbox{\it Wishart}_d\left(N, 0,R\right)$.
By Theorem 2 in~\cite{Steerneman/vanPerlo-tenKleij:2008}, we have
$$
q(z+Z)\in \Sa_d^+\quad \mbox{\rm with probability 1} \Longleftrightarrow
\left(N\geq d\quad \mbox{\rm and} \quad q(z)+R\in \Sa_d^+\right).$$

When $N\geq d$, non-central Wishart matrices can be represented in terms of 
random matrices  
\begin{equation}\label{def-ZZ}
\ZZ:=\left(\ZZ^1,\ldots,\ZZ^N\right)\in \Ma_{d\times N}
\end{equation}
associated with $N$
independent and centered Gaussian column vectors $$\ZZ^i=\left(
\begin{array}{c}
\ZZ^{i}_1\\
\vdots\\
\ZZ^{i}_{d}
\end{array}
\right)\sim \Na_d(0,\I)$$  We recall the multivariate central limit theorem
\begin{equation}\label{def-HH}
 \HH^N:=\frac{1}{\sqrt{N}}
\left( \ZZ\, \ZZ^{\prime}-N~\I\right)\hooklongrightarrow_{N\rightarrow\infty}\sqrt{2}~\GG_{\tiny sym}:=\frac{\GG+\GG^{\prime}}{\sqrt{2}}
\end{equation}
with the $(d\times d)$-submatrix
\begin{equation}\label{def-GG}
\GG:=\left(\ZZ^1,\ldots,\ZZ^d\right)\in \Ma_{d\times d}.
\end{equation}
In Random Matrix Theory, the limiting random matrix $\sqrt{2}\GG_{\tiny sym}$ is sometimes called the Gaussian orthogonal ensemble.

 \begin{lem}\label{lem-1-intro}
 Assume that $N\geq d$ and $q(z),R\in \Sa_d^+$. In this situation,
 we have the orthogonal decomposition
\begin{equation}\label{wish-dec-intro}
 q(z+Z)\stackrel{law}{=}q(z)+R+\frac{1}{\sqrt{N}}~\Delta^N_q(z)
 \end{equation}
 with the remainder term
 $$
 \Delta^N_{q}(z)=R^{1/2}\, \HH^N\,R^{1/2}+2~\left(R^{1/2}~\GG~q^{1/2}(z)~\right)_{\tiny sym}
 $$
 \end{lem}
 The proof of the above lemma is rather technical, thus it is provided in the appendix on page~\pageref{lem-1-intro-proof}. Lemma~\ref{lem-1-intro} can be used to derive several estimates. For instance,  we have the Cauchy-Schwarz inequality
\begin{eqnarray*}
\Vert R^{1/2}\, \HH^N\,R^{1/2}\Vert^2_F&=&\tr\left( \left(\HH^N\,R\, \HH^N\right)\,R)\right) \leq \lambda_{1}(R)^2~\tr( (\HH^N)^2)\\
&=&  \lambda_{1}(R)^2~\sum_{1\leq k,l\leq d}~\left(\frac{1}{\sqrt{N}}~\sum_{1\leq i\leq N} (\ZZ^i_k \ZZ^i_l-1_{k=l})\right)^2
\end{eqnarray*}
In the same vein, we have
\begin{eqnarray*}
\Vert R^{1/2}~ \GG~q^{1/2}(z)\Vert^2_F&\leq& \lambda_{1}(R)\lambda_{1}(q(z))~\tr(\GG\GG^{\prime})
\end{eqnarray*}
This yields for any $r\geq 1$ the estimate
\begin{equation}\label{prop-moments-loc}
\EE\left(\Vert  \Delta^N_{q}(z)\Vert^r_F\right)^{1/r}\leq c(r)~ \left(\lambda_{1}(R)+\sqrt{\lambda_{1}(R)\lambda_{1}(q(z))}\right)
\end{equation}
Working a little harder we also have the variance formulae.
\begin{prop}\label{var-ncW-prop}
For any $N\geq 1$ we have
\begin{equation}\label{var-ncW}
\EE\left( \left( q(z+Z)-(q(z)+R)\right)^2\right)=\frac{1}{N}~ \EE\left(\Delta^N_q(z)^2\right)
\end{equation}
with the matrix
$$
 \EE\left(\Delta^N_q(z)^2\right)=R^{2}+\tr(R)~R+
q(z)R+Rq(z)+\tr(R) q(z)+\tr(q(z))R
$$
\end{prop}
The proof of Proposition~\ref{var-ncW-prop} is provided in the appendix, on page~\pageref{var-ncW-proof}.

\subsection{Orthogonal decompositions}\label{sec-orth-decomp}
Consider the rescaled covariance mapping
$$
 p~:~x\in \Ma_{d\times (N+1)}\mapsto  p(x):=\frac{1}{N}~(x-M(x))(x-M(x))^{\prime}\in \Sa^0_{d}
$$
with the map $x\mapsto M(x)$ defined in (\ref{def-Mmu}).
Let $\Za=(\Za^1,\ldots,\Za^{N+1})\in \RR^{d\times (N+1)}$ be a random matrix with $(N+1)$ independent centered Gaussian  $d$-column vectors
$\Za^i$ with covariance matrix $R>0$; that is for any $1\leq i\leq (N+1)$
$$
\Za^i=\left(
\begin{array}{c}
\Za^{i}_1\\
\vdots\\
\Za^{i}_{d}
\end{array}
\right)\sim \Na_d(0,R)
$$

Also denote by $\ZZ^{0}\in \RR^d$ a centered Gaussian random vectors with unit covariance, independent of matrix $\ZZ$ defined in (\ref{def-ZZ}). We also set
 $$
 \Delta^N(x):=R^{1/2}\, \HH^N\,R^{1/2}+2~\left(R^{1/2}~ \GG~p^{1/2}(x)\right)_{\tiny sym}
 $$
 with $(\HH^N,\GG)$ as in (\ref{def-HH}) and (\ref{def-GG}).

 \begin{lem}\label{lem-2-intro}
 Assume that $ N>d$ and $p(x),R\in \Sa_d^+$. In this situation, up to a change of probability space,
 we have the orthogonal decomposition
\begin{equation}\label{wish-dec-intro-2}
\left\{\begin{array}{rcl}
m(x+\Za)&=&\displaystyle m(x)+\frac{1}{\sqrt{N+1}}~R^{1/2}~\ZZ^0\\
&&\\
 p(x+\Za)&=&\displaystyle p(x)+R+\frac{1}{\sqrt{N}}~\Delta^N(x)\sim \frac{1}{N}~\mbox{\it Wishart}_d\left(N, Np(x),R\right)
 \end{array}\right.
 \end{equation}
 \end{lem}
 The proof of the above lemma is also technical, thus it is provided in the appendix on page~\pageref{lem-2-intro-proof}.

\subsection{Fluctuation analysis}\label{sec-fluct-analysis}
Note that
$$
\HH^N=\frac{1}{\sqrt{N}}
\left( \GG\, \GG^{\prime}-d~\I\right)+\sqrt{1-\frac{d}{N}}~
\frac{1}{\sqrt{N-d}}\left(\sum_{d<i\leq N}\ZZ^i(\ZZ^i)^{\prime}-(N-d)\I\right)
$$
Rewritten in a slightly different form, we have
\begin{equation}\label{def-HH-dec}
\HH^N=\sqrt{\frac{d}{N}}~\HH^d
+\sqrt{1-\frac{d}{N}}~\HH^{(N-d)}_{-d}
 \end{equation}
with the centered random matrix
$$
\HH^{(N-d)}_{-d}:=\frac{1}{\sqrt{N-d}}\left(\sum_{1\leq i\leq (N-d)}\ZZ^{d+i}(\ZZ^{d+i})^{\prime}-(N-d)\I\right)
$$
independent of $\GG=(\ZZ^1,\ldots,\ZZ^d)$ and thus independent of $(\GG,\HH^d)$.
Note that
$$
 \HH^d:=\frac{1}{\sqrt{d}}
\left( \GG\, \GG^{\prime}-d~\I\right)
$$
We also have the decomposition
$$
\HH^N=\HH^{(N-d)}_{-d}+\sqrt{\frac{d}{N}}~\left(\HH^d
+\frac{\left(\sqrt{1-\frac{d}{N}}-1\right)\left(\sqrt{1-\frac{d}{N}}+1\right)}{\sqrt{\frac{d}{N}}\left(\sqrt{1-\frac{d}{N}}+1\right)}~\HH^{(N-d)}_{-d}\right)
$$
This yields the formula
\begin{equation}\label{pert-H}
\HH^N=\HH^{(N-d)}_{-d}+\sqrt{\frac{d}{N}}~\HH^{(N,d)}
\end{equation}
with the random matrix
$$
 \HH^{(N,d)}:=\HH^d
-\frac{1}{\sqrt{\frac{N}{d}-1}+\sqrt{\frac{N}{d}}}~\HH^{(N-d)}_{-d}
$$
The above decompositions yield the following fluctuation result.
\begin{prop} 
We have the multivariate central limit theorem
\begin{equation}\label{def-HH-22}
 (\GG,\HH^N)\quad \mbox{and}\quad  \left(\GG,\HH^{(N-d)}_{-d}\right)
\quad \hooklongrightarrow_{N\rightarrow\infty}\quad(\GG,\HH)
 \end{equation}
 In the above display $\HH$ stands for  an independent copy of the random matrix $\sqrt{2}\GG_{\tiny sym}$.
\end{prop}

\subsection{Inverse matrix moments}\label{sec-inv-moments}

For central Wishart matrices $$q(Z)\sim \frac{1}{N} \mbox{\it Wishart}_d\left(N,R\right)$$ we have the well known (see for instance~\cite{muirhead}) simple formula
\begin{equation}\label{inv-central-W}
\EE(q(Z)^{-1})=\frac{1}{1-(d+1)/N}~R^{-1}
\end{equation}
The one dimensional non central case can also be easily handle in terms of non-central 
chi-squares, see for instance~\cite{horton-2}.  To the best of our knowledge, no explicit closed form expressions nor any inductive formulae are available for multivariate inverse non central Wishart matrices.  
The calculation of expectation of the inverse moment of these random matrices often relies on unresolved infinite series involving zonal polynomials or hypergeometric functions of matrix arguments, see for instance Proposition 4.3 in~\cite{letac-1} and Theorem 1 in~\cite{hillier-kan}. Next we provide a rather elementary way of upperbounding these quantities.

Firstly, note that
$$
\ZZ=(\GG,\ZZ_{-d})\quad\mbox{\rm with}\quad
\ZZ_{-d}:=\left(\ZZ^{d+1},\ldots,\ZZ^{d+(N-d)}\right)
$$
This yields the decomposition
$$
\ZZ\, \ZZ^{\prime}=\GG\, \GG^{\prime}+\ZZ_{-d}\, \ZZ_{-d}^{\prime}
$$
from which we check that
$$
  \begin{array}{l}
\displaystyle q(z+Z)\stackrel{law}{=}q(z)+
 \frac{1}{N}~R^{1/2}
\GG\, \GG^{\prime}R^{1/2} 
 +\frac{1}{\sqrt{N}}~~\left(q^{1/2}(z)~ \GG^{\prime}~R^{1/2}+R^{1/2}~ \GG~q^{1/2}(z)\right) \\
 \\ 
 \hskip3cm\displaystyle
+ \frac{1}{N}~R^{1/2}\ZZ_{-d}\, \ZZ_{-d}^{\prime}R^{1/2}
 \end{array}$$
 Note that
 $$
(N-d)~q(\ZZ_{-d}):= \ZZ_{-d}\, \ZZ_{-d}^{\prime}~\sim~\mbox{\it Wishart}_d\left((N-d),\I\right)
 $$
Rewritten in a slightly different form, we have
 $$
   \begin{array}{l}
\displaystyle
 R^{1/2}\ZZ_{-d}=Z_{-d}:=(Z^{d+1},\ldots,Z^{d+(N-d)})\\
 \\
 \Longrightarrow
 (N-d)~q(Z_{-d})=Z_{-d}\, Z_{-d}^{\prime}~\sim~\mbox{\it Wishart}_d\left((N-d),R\right)
  \end{array}$$
In this notation, we have
$$
  \begin{array}{l}
\displaystyle q(z+Z)\\
\\
\displaystyle\stackrel{law}{=}\left(q(z)^{1/2}+
 \frac{1}{\sqrt{N}}~R^{1/2}~\GG~
\right)\left(q(z)^{1/2}+
 \frac{1}{\sqrt{N}}~R^{1/2}
\GG~
\right)^{\prime}
+ \frac{1}{N}~R^{1/2}\ZZ_{-d}\, \ZZ_{-d}^{\prime}R^{1/2}
 \end{array}$$
We summarize the above discussion with the following lemma.
\begin{lem}
For any  $z\in \Ma_{d\times N}$ we have
 $$
  \begin{array}{l}
\displaystyle q(z+Z)\\
\\
\displaystyle\stackrel{law}{=}\left(q(z)^{1/2}+
 \frac{1}{\sqrt{N}}~R^{1/2}~\GG~
\right)\left(q(z)^{1/2}+
 \frac{1}{\sqrt{N}}~R^{1/2}~
\GG~
\right)^{\prime}
+ \left(1-\frac{d}{N}\right)~ q(Z_{-d})
 \end{array}$$
 \end{lem}
 This yields for any $z\in \Ma_{d\times N}$ the uniform lower bound
 $$
N q(z+Z)\stackrel{law}{\geq} ~ (N-d)~q(Z_{-d})=R^{1/2}~\ZZ_{-d}\, \ZZ_{-d}^{\prime}~R^{1/2}
 $$
Therefore using (\ref{inv-central-W}) for any $N>2d+1$ we have the estimate
 $$
 \EE\left( q(z+Z)^{-1}\right)\leq 
N~R^{-1/2}~\EE\left(\left(\ZZ_{-d}\, \ZZ_{-d}^{\prime}\right)^{-1}\right)~R^{-1/2}=  \frac{N}{(N-d)-(d+1)}~R^{-1}~
 $$
 Recalling that   $P\mapsto P^{-1}$ is convex on $\Sa^+_d$ we check the following lemma.
 \begin{lem}
 For any $N>2d+1$ and $z\in \Ma_{d\times N}$ we have the estimate
 \begin{equation}\label{ref-Inv-ncwishart}
 (q(z)+R)^{-1}\leq  \EE\left( q(z+Z)^{-1}\right)\leq  \frac{R^{-1}}{1-(2d+1)/N}~
 \end{equation}
 \end{lem}

\section{Ensemble covariance estimates}

This section is devoted to a probabilistic and stability analysis of the sample covariance matrices generated by the Ensemble Kalman Filter (EnKF). Our objective is twofold. First, we provide a stochastic perturbation representation of the EnKF covariance recursion, which makes explicit the random fluctuations induced by finite ensemble size. Second, we exploit this representation to establish Lyapunov stability and ergodic properties of the resulting stochastic Riccati dynamics.
The key observation underlying our approach is that the EnKF covariance update–prediction mechanism can be interpreted as a two-step Markov chain on the cone of positive definite matrices, driven by non-central Wishart random matrices. This interpretation allows us to separate the deterministic Riccati flow from the stochastic perturbations arising from sampling noise and to analyze the latter using tools from random matrix theory and Markov process stability.

In the first subsection, we introduce a convenient Gaussian perturbation framework that embeds the EnKF covariance recursion into an equivalent stochastic system. This representation shows that the sample covariance matrices evolve as perturbations of the classical Kalman Riccati equation, with fluctuations of order $N^{-1/2}$. 
In particular, both the updating and prediction steps are expressed in terms of random matrices constructed from independent Gaussian ensembles, leading naturally to non-central Wishart distributions.

Building on this formulation, the second subsection focuses on the long-time behavior of the EnKF covariance process. We prove that the sequence of sample covariance matrices defines a Markov chain that remains almost surely in the cone of positive definite matrices. By introducing a suitable Lyapunov function that controls both the trace and the inverse trace of the covariance matrices, we establish geometric ergodicity of the chain under mild conditions on the ensemble size. As a consequence, the stochastic Riccati dynamics admits a unique invariant distribution, and convergence toward equilibrium holds exponentially fast in weighted total variation norms.

These results provide a rigorous probabilistic foundation for the stability of the EnKF covariance estimates and clarify the role played by ensemble size in controlling stochastic fluctuations. They also offer a useful framework for further quantitative error analysis and for the study of more general ensemble-based filtering algorithms.

\subsection{A stochastic perturbation theorem}

Consider a collection of independent
random matrices  $$
\ZZ_n:=\left(\ZZ^{1}_n,\ldots,\ZZ^{N}_n\right)\in \Ma_{d\times N}
\quad\mbox{\rm and}\quad
\widehat{\ZZ}_n:=\left(\widehat{\ZZ}^{1}_n,\ldots,\widehat{\ZZ}^{N}_n\right)\in \Ma_{d_0\times N}
$$
indexed by $n\in\NN$ associated with $N$
independent and centered Gaussian column vectors $\ZZ^{i}_n\sim \Na_d(0,\I)$ and
$\widehat{\ZZ}^{i}_n\sim \Na_{d_0}(0,\I)$. We also denote by
$(\GG_n,\HH^N_{n})$ and respectively $(\widehat{\GG}_n,\widehat{\HH}^N_{n})$, the random matrices defined in  (\ref{def-HH}) and (\ref{def-GG}) 
by replacing $\ZZ$ by $\ZZ_n$ and respectively by $\widehat{\ZZ}_n$. Also denote by $$(\ZZ^{0}_n,\widehat{\ZZ}^{0}_n)\in (\RR^d\times \RR^{d_0})$$ a sequence of independent centered Gaussian random vectors with unit covariance, independent of the sequence $(\ZZ_n,\widehat{\ZZ}_n)_{n\geq 0}$.
In this notation we have the  lemma.
\begin{lem}
Up to a change of probability space we have
 \begin{eqnarray*}
m_0&=&\EE(X_0)+\frac{1}{\sqrt{N+1}}~P_0^{1/2}\,\ZZ^{0}_0
\\
p_0&=&P_0+\frac{1}{\sqrt{N}}~ \Lambda_{0}^N
\quad \mbox{with}\quad
 \Lambda_{0}^N:=P_0^{1/2}\, \HH^N_{0}\,P_0^{1/2}
\end{eqnarray*}
\end{lem}
For any $n\geq 0$, we  set
 \begin{eqnarray*}
 \Delta^N_{n+1}&:=&R^{1/2}\, \HH^N_{n+1}\,R^{1/2}+2~\left(R^{1/2}~ \GG_{n+1}~(A\,\widehat{p}_{n}~A^{\prime})^{1/2}\right)_{\tiny sym}
\\
 \widehat{\Delta}^N_{n}&:=&\widehat{R}(p_n)^{1/2}~\widehat{\HH}^N_{n}~\widehat{R}(p_n)^{1/2}+2~\left(\widehat{R}(p_n)^{1/2}
~\widehat{\GG}_n~\left(
\widehat{A}(p_n)~p_n~\widehat{A}(p_n)^{\prime}\right)^{1/2}\right)_{\tiny sym}
\end{eqnarray*}
with the mappings $(\widehat{A}(p),\widehat{R}(p))$ defined in (\ref{eq:Gn}) and (\ref{eq:gnB-ident1}). By (\ref{matrix-form}), the updating and the prediction of the EnKF take the following matrix form
\begin{eqnarray*}
\widehat{\xi}_{n}
&=&\widehat{A}(p_n)~\xi_n +K(p_n)~\YY_n-K(p_n)\VV_n
\\
\xi_{n+1}&=&A ~\widehat{\xi}_{n}+ \WW_{n+1}\end{eqnarray*}
Recall from (\ref{eq:gnB-ident1})
that for any $p\in\Sa_d^0$ we have the decomposition
$$
\widehat{A}(p)p=\widehat{A}(p)~p~\widehat{A}(p)^{\prime}+\widehat{R}(p)
\quad \mbox{\rm and}\quad  \widehat{R}(p):=K(p)\,R_0\,K(p)^{\prime}
$$

Next theorem is now a  consequence of Lemma~\ref{lem-2-intro}. 
\begin{theo}\label{theo-mp-intro}
Up to a change of probability space, the updating transition  of the EnKF is given for any $n\geq 0$ by the equations
$$
 \left\{
\begin{array}{rcl}
\widehat{m}_{n}&=&\displaystyle m_n+K(p_n)
~(Y_n-B\,  m_n) +\frac{1}{\sqrt{N+1}}~K(p_n)~R_0^{1/2}~\widehat{\ZZ}^{0}_n\\
 \widehat{p}_n&=&\displaystyle\widehat{A}(p_n)~ p_n+\frac{1}{\sqrt{N}}~ \widehat{\Delta}^N_{n}
\end{array}\right.
$$
Moreover, the prediction step of the EnKF is given  for any $n\geq 0$ by the equations
$$
 \left\{
\begin{array}{rcl}
m_{n+1}&=&\displaystyle A ~\widehat{m}_{n}+\frac{1}{\sqrt{N+1}}~R^{1/2}~\ZZ^{0}_{n+1}\\
 p_{n+1}
&=&\displaystyle A ~\widehat{p}_{n}~A ^{\prime}+R+\frac{1}{\sqrt{N}}~\Delta^N_{n+1}
\end{array}\right.
$$
\end{theo}
This yields the following corollary.
\begin{cor}\label{cor:wishart}
For any $n\geq 0$ we have
\begin{equation}\label{eq:pn}
 p_{n}
=\Phi(p_{n-1})+\frac{1}{\sqrt{N}}~\Lambda^N_{n}\quad \mbox{with}\quad \Lambda^N_{n}:=
 A~\widehat{\Delta}^N_{n-1}~A ^{\prime}+\Delta^N_{n}.
\end{equation}
with the convention $\Phi(p_{-1})=P_0$ for $n=0$
as well as
$$ \widehat{\Delta}^N_{-1}=0\quad
\mbox{ and}\quad
\Delta^N_{0}:=P_0^{1/2}\, \HH^N_{0}\,P_0^{1/2}.$$
\end{cor}

\subsection{Lyapunov stability analysis}

The next Theorem is a consequence of Theorem~\ref{theo-mp-intro}. It provides an interpretation of the updating-prediction steps of the EnKF sample covariances  in terms of a two-steps Markov chain involving non-central Wishart  random matrices and the mappings $P\mapsto \widehat{A}(P)$ and $P\mapsto\widehat{R}(P)$ defined in (\ref{eq:Gn}) and (\ref{eq:gnB-ident1}).
To avoid unnecessary technical discussions, without further mention, we assume that $P_0,R_0,R>0$ and $N+1>d$.
\begin{theo}\label{stoch-ricc-theo}
The sample covariances $p_n, \widehat{p}_n$ of the EnKF belongs to $\Sa^+_d$ with probability $1$ and they are  given by the initial condition
$$
p_0\sim \frac{1}{N}~ \mbox{\it Wishart}_{d}\left(N, 0,P_0 \right)
$$
and for any $n\geq 0$ by the updating-prediction steps
$$
\left\{
\begin{array}{rcl}
 \widehat{p}_n&\sim&\displaystyle \frac{1}{N}~ \mbox{\it Wishart}_{d}\left(N, N~(\widehat{A}(p_n)~p_{n}~
\widehat{A}(p_n)^{\prime}),\widehat{R}(p_n) \right)\\
&&\\
 p_{n+1}&\sim &\displaystyle\frac{1}{N}~ \mbox{\it Wishart}_{d}\left(N, N(A ~\widehat{p}_{n}~A ^{\prime}),R\right)
\end{array}\right.$$
\end{theo}

Given some non-negative function $\Ua\geq 1$ on $\Sa_d^0$ we define the $\Ua$-norm of some function $f: \Sa_d^0 \to \RR$ and some locally finite signed measure $\mu(dp)$ on $\Sa_d^0$  by the formulae
\begin{equation}\label{U-norm}
\Vert f\Vert_{\tiny \Ua}:=\sup_{p\in \Sa_d^0}\left\vert {f(p)}/{\Ua(p)}\right\vert\quad \mbox{\rm and}\quad \Vert \mu\Vert_{\tiny \Ua}:=\sup{\{\vert \mu(f)\vert~:~f~\mbox{\rm s.t.}~\Vert f\Vert_{\tiny \Ua}\leq 1\}}.
\end{equation}
Let $\Pa_{\Ua}(\Sa_d^0)$ be the set of probability measures $\mu$ on $\Sa_d^0$ such that $\mu(\Ua)<\infty$.

Let $\Ta(p, dq)$ denote the Markov transition of the sample covariance matrices
  $p_n$ presented in Theorem~\ref{stoch-ricc-theo}; that is,
$$
\Ta(p, dq)=\PP\left(p_{n+1}\in dq~|~p_{n}=p\right)
$$
We are using $\Ta$ to denote both the left-action operator on measure
 $\mu\mapsto \mu\Ta$ and the right-action operator on functions
$f\mapsto\Ta(f)$  defined by
$$
(\mu\Ta)(dq)=\int~\mu(dp)~\Ta(p,dq)\quad\mbox{\rm and}\quad
\Ta(f)(p):=\int\Ta(p,dq) f(q)
$$
The corresponding transition semigroup $(\Ta_n)_{n\geq 0}$ is defined by
$$
\Ta_n(p, dq)=\PP\left(p_{n}\in dq~|~p_{0}=p\right)
$$

Equivalently, the integral operators
$\Ta_{n+1}$ are defined
for any $n\geq 0$ by the induction
$$
\Ta_{n+1}=\Ta\circ\Ta_{n}=\Ta_{n}\circ \Ta
$$
with the convention $\Ta_0(f)=f$ and $\mu\Ta_0=\mu$, for $n=0$.
In the above display, the symbol $``\circ"$ stands for the composition of integral operators.

Consider the function with compact level subsets $\{\Ua\leq l\}$ on $\Sa_d^0$ defined by
$$
\Ua(p):=1+\tr(p)+\tr(p^{-1})
$$
In proposition~\ref{prop-estimates-pi} we prove that
$
\Vert\Ta(\Ua)\Vert<\infty
$, as soon as $N\geq (2d+1)$. This property ensures that $\Ua$ is a Lyapunov function of $\Ta$ in the sense there exists some  $\epsilon\in ]0,1[$ and some constant $c>0$ such that
$$
\Ta(\Ua)\leq \epsilon~\Ua+c
$$
Applying Theorem 2.2 and Remark 2.3 in~\cite{marc+elmaati} we obtain the following stability theorem for the stochastic Riccati equations presented in Theorem~\ref{stoch-ricc-theo}.
\begin{theo}
Assume that $N\geq (2d+1)$. In this situation, there exists an unique invariant measure $\pi=\pi\Ta\in \Pa_{\Ua}(\Sa_d^0)$.
In addition, there exists some parameters $a,b>0$ such that for any $n\geq 1$ and any $\mu_1,\mu_2\in \Pa_{\Ua}(\Sa_d^0)$ we have the $\Ua$-norm contraction estimates
$$
\Vert \mu_1\Ta_n- \mu_2\Ta_n \Vert_{\Ua}\leq a~e^{-bn}~\Vert \mu_1- \mu_2\Vert_{\Ua}
$$
as well as  the exponential decays w.r.t. the total variation norm
$$
\Vert \mu_1\Ta_n- \mu_2\Ta_n \Vert_{\tiny tv}\leq a~e^{-bn}~\Vert \mu_1- \mu_2\Vert_{\tiny tv}.
$$
\end{theo}

\section{Stochastic perturbation analysis }

This section analyzes the stochastic fluctuations of the Ensemble Kalman Filter (EnKF) arising from finite ensemble size. We first establish time-uniform bounds for the sample covariances, showing that they remain positive definite and their moments are uniformly controlled over time. These results ensure stability and provide a foundation for bias and fluctuation analysis.
We then examine the bias of the sample covariance, demonstrating that its expected value is systematically below the deterministic Riccati solution. Using second-order expansions, we isolate leading-order stochastic perturbations and higher-order bias terms, explicitly quantifying the effect of ensemble size.
Finally, we derive central limit theorems for the covariance fluctuations and time-uniform estimates for the state errors, showing that both converge with order $1/\sqrt{N}$. These results provide a rigorous probabilistic description of the EnKF, capturing both random variability and finite-ensemble bias, and support precise error control in practical implementations.

\subsection{Some time uniform estimates}

\begin{prop}\label{prop-estimates-pi}
Assume that $P_0,R_0,R>0$ and $N+1>d$. In this situation, we have $p_n>0$ and $\widehat{p}_n>0$. For any $r\geq 1$, there also exist some finite constants $c_1(r),c_2(r)$ such that for any $n\geq 0$ we have the almost sure time uniform estimates
\begin{equation}\label{unif-est}
\EE\left(\Vert\widehat{\Delta}^N_n\Vert^{r}~|~\xi_n\right)\vee
\EE\left(\Vert\Delta^N_{n+1} \Vert^{r}~|~\xi_n\right)\vee \EE\left(\Vert\Lambda^N_{n+1}\Vert^{r}~|~\xi_n\right)\leq c_1(r)
\end{equation}
as well as
$$
\EE\left(\left( \tr( \widehat{p}_n) \right)^{r}~|~\xi_n\right)\vee \EE\left(\left( \tr( p_{n+1}) \right)^{r}~|~\xi_n\right)\leq c_2(r)
$$
In addition, for any $N\geq (2d+1)$ and $n\geq 1$ we have
\begin{equation}\label{unif-est-2}
R\leq  \EE(p_{n})\leq AS^{-1}A^{\prime}+R\quad \mbox{and}\quad 
  \EE(p_{n}^{-1})\leq  \frac{R^{-1}}{1-(2d+1)/N}
\end{equation}
\end{prop}

\proof

The article \cite{Dykstra:1970} proved that $p_0$ is non-singular if $N+1>d$ and $P_0>0$. Recalling that given $\widehat{p}_{n}$, the sample covariance
$p_{n+1}$ has noncentral Wishart distribution, by Theorem 2 in~\cite{Steerneman/vanPerlo-tenKleij:2008} whenever $N+1>d$ we also check that
$$
R >0\Longrightarrow \forall n\geq 0\quad
A ~\widehat{p}_{n}~A ^{\prime}+R >0\Longrightarrow   \forall n\geq 1\quad p_{n}>0
$$
In addition, we have
$$
\begin{array}{l}
p_{n}>0
 \Longrightarrow
(\I-K(p_n)B)p_n=(p_n^{-1}+S)^{-1}>0
\Longrightarrow \widehat{p}_{n}>0
\end{array}
$$
 This ends the proof of the first assertion.
 Observe that
 $$
R\leq  \EE(p_{n+1})=\EE(\Phi(p_n))\leq AS^{-1}A^{\prime}+R
 $$
On the other hand, using (\ref{ref-Inv-ncwishart})  we also have
$$
 \frac{R^{-1}}{1-(2d+1)/N}~
\geq \EE(p_{n+1}^{-1})\geq (\EE(p_{n+1}))^{-1}\geq (AS^{-1}A^{\prime}+R)^{-1}
$$
Now we come to the proof of \eqref{unif-est}. Note that
$$
\widehat{R}(P):=K(P)\,R_0\,K(P)^{\prime}\leq \lambda_{1}(R_0)~K(P)K(P)^{\prime}
$$
For any $P>0$ we have
\begin{eqnarray*}
K(P)B &=&\I-(\I+PS)^{-1}=(\I+PS)^{-1}PS=(P^{-1}+S)^{-1}S
\end{eqnarray*}
This yields the estimate
$$
\lambda_{d}(BB^{\prime})~K(P)K(P)^{\prime}\leq K(P)BB^{\prime}K(P)^{\prime}\leq \left({\lambda_{1}(S)}/{\lambda_{d} (S)}\right)^2~\I
$$
from which we check that
$$
\lambda_{1}(\widehat{R}(P))\leq \frac{\lambda_{1}(R_0)}{\lambda_{d}(BB^{\prime})~}~\left(\frac{\lambda_{1}(S)}{\lambda_{d} (S)}\right)^2
$$
In the same vein, we have
\begin{eqnarray*}
\widehat{A}(P)P\widehat{A}(P)^{\prime}
&=&(\I+PS)^{-1}P(\I+SP)^{-1}=(P^{-1}+S)^{-1}(S^{-1}+P)^{-1}S^{-1}
\end{eqnarray*}
and therefore
\begin{eqnarray*}
\tr(\widehat{A}(P)P\widehat{A}(P)^{\prime})&\leq& \frac{1}{\lambda_{d}(S)}~\tr((P^{-1}+S)^{-1}(S^{-1}+P)^{-1})\\
&\leq &\frac{1}{\lambda_{d}(S)^2}~\tr((S^{-1}+P)^{-1})\leq \frac{1}{\lambda_{d}(S)^2}~\tr(S)
\end{eqnarray*}

By (\ref{prop-moments-loc})  we have the time-uniform almost sure estimate
$$
\begin{array}{l}
\EE\left(\Vert \widehat{\Delta}^N_n\Vert_F^r~|~p_n\right)^{1/r}
\leq c_1(r)
\end{array}$$
In the same vein, we have
$$
\begin{array}{l}
\EE\left(\Vert  \Delta^N_{n+1}\Vert_F^r~|~\widehat{p}_n\right)^{1/r}
\leq c_3(r)~\left(1+\sqrt{\tr(\widehat{p}_{n})}\right)
\end{array}$$
On the other hand, we have
$$
\tr( \widehat{p}_n)=\tr((\I+p_nS)^{-1} p_n)+\frac{1}{\sqrt{N}}~ \tr(\widehat{\Delta}^N_n
)\leq \tr(S^{-1})+\frac{1}{\sqrt{N}}~ \tr(\widehat{\Delta}^N_n
)
$$
and by a direct application of Cauchy-Schwarz inequality we have
$$
\vert \tr(\widehat{\Delta}^N_n)\vert\leq c_4~\Vert\widehat{\Delta}^N_n\Vert_F
\quad\mbox{\rm
and therefore}\quad
\EE\left(\Vert   \Lambda^N_{n+1}\Vert_F^r~|~p_n\right)^{1/r}\leq c_6(r)
$$
\qed
\subsection{Under bias property}
Recalling that the mapping $$p\mapsto (\I+p^{-1})^{-1}$$ is concave on the set of positive matrices (cf. for instance Theorem 4.1.1 in~\cite{bathia}) we check that
\begin{eqnarray*}
\EE\left(A ~(\I+p_n S)^{-1}~ p_n~A ^{\prime}\right)&=&A~S^{-1/2} \EE\left((S^{-1/2}p_n^{-1}S^{-1/2}+\I)^{-1}~ \right)~S^{-1/2} A^{\prime}\\
&\leq &A~S^{-1/2} (S^{-1/2}\EE\left(p_n\right)^{-1}S^{-1/2}+\I)^{-1}~~S^{-1/2} A ^{\prime}
\end{eqnarray*}
from which we check that for any $n\geq 0$ we have $$
R\leq \EE\left(p_{n+1}\right) =\EE\left(\Phi(p_n)\right)  \leq \Phi(\EE(p_n))\leq R+A~S^{-1} A^{\prime}.
$$
Using Lemma~\ref{lem-phi-props} this implies that
$$
\EE\left(p_{n+2}\right) =\EE\left(\Phi(p_{n+1})\right)  \leq \Phi(\EE(p_{n+1}))\leq
\Phi_2(\EE(p_n))
$$
Iterating the above procedure and recalling that $\EE(p_0)=P_0$ we conclude
that $p_n$ is under biased, that is for any $n\geq 1$ we have
$$
R\leq \EE\left(p_{n}\right)\leq P_n\leq AS^{-1}A^{\prime}+R.
$$

\subsection{Second order expansions}\label{th-intro-punif-proof}
This section is mainly concerned with the proof of Theorem~\ref{th-intro-punif}.

Using the telescoping decomposition
$$
p_n-\Phi_n(P_0)=\Phi_n(p_0)-\Phi_n(P_0)+\sum_{1\leq k\leq n}
\left(\Phi_{n-k}\left(p_k\right)-\Phi_{n-(k-1)}\left(p_{k-1}\right)\right)
$$
we check that
$$
\begin{array}{l}
\displaystyle\ p_n-P_n
=\Phi_n\left(P_0+\frac{1}{\sqrt{N}}~ \Lambda^N_{0}\right)-\Phi_n(P_0)\\
\\
\displaystyle\hskip3cm +\sum_{1\leq k\leq n}
\left(\Phi_{n-k}\left(\Phi(p_{k-1})+\frac{1}{\sqrt{N}}~ \Lambda^N_{k}\right)-\Phi_{n-k}\left(\Phi(p_{k-1})\right)\right)
\end{array}
$$
Using \eqref{phi-form},  we find that
$$
\Phi_n\left(P_0+\frac{1}{\sqrt{N}}~ \Lambda^N_0\right)-\Phi_n(P_0)=\frac{1}{\sqrt{N}}~
\Ea_n(p_0)~ \Lambda^N_0~\Ea_n(P_0)^{\prime}
$$
as well as
$$
\begin{array}{l}
\displaystyle
\Phi_{n-k}\left(\Phi(p_{k-1})+\frac{1}{\sqrt{N}}~ \Lambda^N_{k}\right)-\Phi_{n-k}\left(\Phi(p_{k-1})\right)=\frac{1}{\sqrt{N}}~\Ea_{n-k}(p_k)~ \Lambda^N_{k}~\Ea_{n-k}(\Phi(p_{k-1}))^{\prime}
\end{array}
$$
This yields the following lemma.
\begin{lem}
For any $n\geq 0$ we have interpolating formula
$$
\begin{array}{l}
\displaystyle\ p_n=P_n
+\frac{1}{\sqrt{N}}~ \sum_{0\leq k\leq n}
\Ea_{n-k}(p_k)~ \Lambda^N_{k}~\Ea_{n-k}(\Phi(p_{k-1}))^{\prime}
\end{array}
$$
with the convention $\Phi(p_{-1})=P_0$ for $k=0$. 
\end{lem}
Using the Floquet type formula \eqref{form-Floquet}, for any sub-multiplicative matrix norm (such as for the spectral norm or the Frobenius norm) we have
\begin{equation}\label{unif-pP}
\Vert \ p_n-P_n\Vert\leq \frac{\iota^2}{\sqrt{N}}~ \sum_{0\leq k\leq n}
\Vert\Ea(P_{\infty})\Vert^{2(n-k)}~\Vert \Lambda^N_{k}\Vert
\end{equation}
Using  \eqref{unif-est}, we readily check (\ref{unif-cov}).

Using Proposition~\ref{prop-dec-Ea} we have the first order expansion
$$
\Ea_{n-k}(p_k)=\Ea_{n-k}(\Phi(p_{k-1}))+\nabla \Ea_{n-k}(p_k,\Phi(p_{k-1}))\cdot (p_k-\Phi(p_{k-1}))
$$
This yields the following proposition.
\begin{prop}
For any $n\geq 0$ we have the second order decomposition
\begin{equation}\label{2nd-order}
\begin{array}{l}
\displaystyle\ p_n=P_n+\frac{1}{\sqrt{N}}~ \sum_{0\leq k\leq n}
\Ea_{n-k}(\Phi(p_{k-1}))~ \Lambda^N_{k}~\Ea_{n-k}(\Phi(p_{k-1}))^{\prime}\\
\\
\hskip3cm\displaystyle-\frac{1}{N}~ \sum_{0\leq k\leq n}~
\Ea_{n-k}(p_k)~ \Lambda^N_{k}~~\Ga_{n-k}~
\La_{n-k}(\Phi(p_{k-1}))^{-1}~ \Lambda^N_{k}~\Ea_{n-k}(\Phi(p_{k-1}))^{\prime}
\end{array}
\end{equation}
\end{prop}
Using (\ref{ref-Lan}) we check that
$$
\Ga_{n-k}~
\La_{n-k}(\Phi(p_{k-1}))^{-1}=\left(\Phi(p_{k-1})+
\left(\Ga_{n-k}^{-1}-P_{\infty}\right)\right)^{-1}~
$$
Taking the expectations in the above display, we find that
\begin{equation}\label{2nd-order-bias}
\begin{array}{l}
\displaystyle\ P_n=\EE(p_n)+\frac{1}{N}~ \sum_{0\leq k\leq n}~
\EE\left(\Ea_{n-k}(p_k)~ \Lambda^N_{k}~~\Ga_{n-k}~
\La_{n-k}(\Phi(p_{k-1}))^{-1}~ \Lambda^N_{k}~\Ea_{n-k}(\Phi(p_{k-1}))^{\prime}\right)
\end{array}
\end{equation}

Arguing as in the proof of (\ref{unif-pP}), there exists some parameter $\iota>0$ such that for any $n\geq 0$ we have
$$
0\leq P_n-\EE(p_n)\leq \frac{\iota}{N}~\I
$$

\subsection{Some fluctuation theorems}
Let $((\widehat{\GG}_k,\widehat{\HH}_k),(\GG_{k},\HH_k))_{k\geq 0}$ be independent copies of the random matrices $(\GG,\HH)$ discussed in (\ref{def-HH-22}). Consider the collection of random matrices
$$
\Lambda_{n}:=
 A~\widehat{\Gamma}_{n-1}~A ^{\prime}+\Gamma_{n}
 $$
 with
 \begin{eqnarray*}
\Gamma_{n}&:=&R_n^{1/2}\, \HH_{n}\,R_n^{1/2}+2~\left(R^{1/2}~ \GG_{n}~(A\,\widehat{P}_{n-1}~A^{\prime})^{1/2}\right)_{\tiny sym}
\\
 \widehat{\Gamma}_{n}&:=&\widehat{R}(P_n)^{1/2}~\widehat{\HH}_{n}~\widehat{R}(P_n)^{1/2}+2~\left(\widehat{R}(P_n)^{1/2}
~\widehat{\GG}_n~\left(
\widehat{A}(P_n)~P_n~\widehat{A}(P_n)^{\prime}\right)^{1/2}\right)_{\tiny sym}
\end{eqnarray*}
with $R_n=1_{n\geq 1}~R+1_{n=0}P_0$. For $n=0$ we also use the convention $\widehat{P}_{-1}=0$. 
\begin{theo}\label{theo-TCL-riccati}
For any time horizon $n\geq 0$ we have the weak convergence
$$
\sqrt{N}\left( p_n-P_n\right)\quad\hooklongrightarrow_{N\rightarrow\infty}\quad
 \sum_{0\leq k\leq n}
\Ea_{n-k}(P_k)~ \Lambda_{k}~\Ea_{n-k}(P_k)^{\prime}
$$
as well as the limiting bias formula
$$
\begin{array}{l}
\displaystyle \lim_{N\rightarrow\infty}N(P_n-\EE(p_n)) \\
\\
\displaystyle=\sum_{0\leq k\leq n}~
\EE\left(\Ea_{n-k}(P_k)~ \Lambda_{k}~\left(P_k+
\left(\Ga_{n-k}^{-1}-P_{\infty}\right)\right)^{-1}~ \Lambda_{k}~\Ea_{n-k}(P_k)^{\prime}\right)>0
\end{array}
$$
\end{theo}

The proof of the above theorem is based on the continuous mapping theorem combined with the second order Taylor expansions presented in (\ref{2nd-order}) and (\ref{2nd-order-bias}). Before entering into the details of the proof note that the central limit theorem (\ref{def-HH-22}) yields for any time horizon $n\geq 0$ the weak convergence
$$
\left((\GG_k,\HH^{(N-d)}_{-d,k}),(\widehat{\GG}_k,\widehat{\HH}^{(N-d)}_{-d,k})\right)_{0\leq k\leq n}\hooklongrightarrow_{N\rightarrow\infty}\quad((\GG_{k},\HH_{k}),(\widehat{\GG}_k,\widehat{\HH}_k))_{0\leq k\leq n}
$$

Denote by
$(\HH^d_{n},\HH^{(N-d)}_{-d,n},\HH^{(N,d)}_{n})$ and respectively $(\widehat{\HH}^d_{n},\widehat{\HH}^{(N-d)}_{-d,n},\widehat{\HH}^{(N,d)}_{n})$, the random matrices defined as in (\ref{def-HH-dec}) and (\ref{pert-H})
by replacing $\ZZ$ by $\ZZ_n$ and respectively by $\widehat{\ZZ}_n$.

We also consider the collection of independent random matrices
 \begin{eqnarray*}
\Gamma^N_{n}&:=&R_n^{1/2}\, \HH^{N-d}_{-d,n}\,R_n^{1/2}+2~\left(R^{1/2}~ \GG_{n}~(A\,\widehat{P}_{n-1}~A^{\prime})^{1/2}\right)_{\tiny sym}
\\
 \widehat{\Gamma}^N_{n}&:=&\widehat{R}(P_n)^{1/2}~\widehat{\HH}^{(N-d)}_{-d,n}~\widehat{R}(P_n)^{1/2}+2~\left(\widehat{R}(P_n)^{1/2}
~\widehat{\GG}_n~\left(
\widehat{A}(P_n)~P_n~\widehat{A}(P_n)^{\prime}\right)^{1/2}\right)_{\tiny sym}
\end{eqnarray*}
with $R_n=1_{n\geq 1}~R+1_{n=0}P_0$, for $n=0$ we also use the convention $\widehat{P}_{-1}=0$. For any time horizon $n\geq 0$ we have the weak convergence
$$
\left(\Gamma^N_{k}, \widehat{\Gamma}^N_{k}\right)_{0\leq k\leq n}\quad\hooklongrightarrow_{N\rightarrow\infty}\quad(\Gamma_{k}, \widehat{\Gamma}_{k})_{0\leq k\leq n}
$$
Combining  the continuous mapping theorem with the second order expansions (\ref{2nd-order}) and (\ref{2nd-order-bias}), Theorem~\ref{theo-TCL-riccati} is  now a direct consequence of  the following theorem.
\begin{theo}
For any time horizon $n\geq 0$ we have the weak convergence
$$
\left(\Delta^N_{k}, \widehat{\Delta}^N_{k}\right)_{0\leq k\leq n}\hooklongrightarrow_{N\rightarrow\infty}\quad(\Gamma_{k}, \widehat{\Gamma}_{k})_{0\leq k\leq n}
\quad\mbox{and}
\quad
\left(\Lambda_{k}^N\right)_{0\leq k\leq n}\hooklongrightarrow_{N\rightarrow\infty}
\left(\Lambda_{k}\right)_{0\leq k\leq n}
$$
\end{theo}
The above theorem is itself a direct consequence of the following lemma combined with  the continuous
mapping theorem.
\begin{lem}
For any $r\geq 1$, any $P_0>0$ and $N>(1+d)$  we have  the time-uniform estimates
$$
\sup_{n\geq 1}\left(\EE\left(\left\Vert\,\Delta^N_{n}- \Gamma^N_{n}\right\Vert^r\right)^{1/r}\vee \EE\left(\left\Vert\,\widehat{\Delta}^N_{n}- \widehat{\Gamma}^N_{n}\right\Vert^r\right)^{1/r}\right)\leq c(r)/\sqrt{N}
$$

\end{lem}
\proof
Following the proof of Proposition~\ref{prop-estimates-pi}, for any $P>0$ we have the uniform estimates
$$
\tr(\widehat{A}(P)P\widehat{A}(P)^{\prime})\leq \frac{\tr(S)}{\lambda_{d}(S)^2}\quad \mbox{\rm and}\quad \tr(\widehat{R}(P))\leq d~\frac{\lambda_{d_0}(R_0)}{\lambda_{d}(BB^{\prime})}~\left({\lambda_{1}(S)}/{\lambda_{d} (S)}\right)^2
$$
On the other hand, for any $P\geq R$ we have $(P^{-1}+S)^{-1}\geq (R^{-1}+S)^{-1}$ and therefore
\begin{eqnarray*}
 K(P)BB^{\prime}K(P)^{\prime}&=&(P^{-1}+S)^{-1}S^2(P^{-1}+S)^{-1}\\
 &\geq &\lambda_d(S)^2~
 (P^{-1}+S)^{-1/2} (P^{-1}+S)^{-1} (P^{-1}+S)^{-1/2}\\
 &\geq &\frac{\lambda_d(S)^2}{\lambda_1(R^{-1}+S)}~ (R^{-1}+S)^{-1}
\end{eqnarray*}
This implies that
\begin{eqnarray*}
P\geq R&\Longrightarrow& \widehat{R}(P)\geq  \frac{\lambda_{d_0}(R_0)\lambda_d(S)^2}{\lambda_1(BB^{\prime})\lambda_1(R^{-1}+S)}~(R^{-1}+S)^{-1}\\
&\Longrightarrow& \widehat{R}(P)^{1/2}\geq \left( \frac{\lambda_{d_0}(R_0)\lambda_d(S)^2}{\lambda_1(BB^{\prime})\lambda_1(R^{-1}+S)}\right)^{1/2}~(R^{-1}+S)^{-1/2}>0
\end{eqnarray*}
For any $R\leq P\leq A ~S^{-1}~A^{\prime}+R$, we also have
\begin{eqnarray*}
\widehat{A}(P)~P~\widehat{A}(P)^{\prime}&=&(\I+PS)^{-1}P(\I+SP)^{-1}\\
&=& (P^{-1}+S)^{-1}P^{-1}(P^{-1}+S)^{-1}\\
&\geq&\lambda_d(A ~S^{-1}~A ^{\prime}+R)~(P^{-1}+S)^{-1}(P^{-1}+S)^{-1}
\end{eqnarray*}
Arguing as above we conclude that
$$
\widehat{A}(P)~P~\widehat{A}(P)^{\prime}\geq 
\frac{\lambda_d(A ~S^{-1}~A ^{\prime}+R)}{\lambda_1(R^{-1}+S)}~(R^{-1}+S)^{-1}
$$
In summary, we have proved that
$$
\begin{array}{l}
R\leq P\leq A ~S^{-1}~A^{\prime}+R\\
\\
\displaystyle\Longrightarrow
\left(\widehat{A}(P)~P~\widehat{A}(P)^{\prime} \right)^{1/2}\geq 
\frac{\lambda_d(A ~S^{-1}~A ^{\prime}+R)}{\lambda_1(R^{-1}+S)}~(R^{-1}+S)^{-1/2}>0
\end{array}$$
Using the Ando-Hemmen inequality
(\ref{square-root-key-estimate}) we check that
$$
\Vert \left(\widehat{A}(p_n)~p_n~\widehat{A}(p_n)^{\prime} \right)^{1/2}- \left(\widehat{A}(P_n)~P_n~\widehat{A}(P_n)^{\prime} \right)^{1/2}\Vert \leq c_1~\Vert p_n- P_n\Vert
$$
as well as
$$
\Vert \widehat{R}(p_n)^{1/2}- \widehat{R}(P_n)^{1/2}\Vert \leq c_2~\Vert p_n- P_n\Vert
$$
Using (\ref{unif-cov}), for any $r\geq 1$ and $N>(1+d)$  we also check  the time-uniform estimates
$$
\sup_{n\geq 1}\EE\left(\left\Vert\, \left(\widehat{A}(p_n)~p_n~\widehat{A}(p_n)^{\prime} \right)^{1/2}- \left(\widehat{A}(P_n)~P_n~\widehat{A}(P_n)^{\prime} \right)^{1/2}\right\Vert^r\right)^{1/r}\leq c_1(r)/\sqrt{N}
$$
as well as
$$
\sup_{n\geq 1}\EE\left(\left\Vert\, \widehat{R}(p_n)^{1/2}- \widehat{R}(P_n)^{1/2}\right\Vert^r\right)^{1/r}\leq c_2(r)/\sqrt{N}
$$
Observe that for any $n\geq 0$ we have
$$
\begin{array}{l}
\Delta_{n}^N- \Gamma^N_{n}\\
 \\
 \displaystyle=\sqrt{\frac{d}{N}}~R_n^{1/2}\, \HH^{(N,d)}_{n}\,R_n^{1/2}+2~\left(R^{1/2}~ \GG_{n}~
 \left((A\,\widehat{p}_{n-1}~A^{\prime})^{1/2}-(A\,\widehat{P}_{n-1}~A^{\prime})^{1/2}\right)\right)_{\tiny sym}
\end{array}$$
with the convention $\widehat{p}_{-1}=0=\widehat{P}_{-1}$ for $n=0$.
The above estimates ensure that
$$
\sup_{n\geq 0}\EE\left(\left\Vert\,\Delta^N_{n}- \Gamma^N_{n}\right\Vert^r\right)^{1/r}\leq c(r)/\sqrt{N}
$$
The end of the lemma now follows elementary computations thus it is skipped.\cqfd

\subsection{State estimates}

\begin{theo}
Assume that $\Vert S^{1/2}AS^{-1/2}\Vert_2<1$. In this situation, for any $r\geq 1$ and $N+1>d$ we have the time-uniform estimates
\begin{equation}\label{stab-unif-theo}
\sup_{n\geq 0}\left(\EE\left(\Vert\, m_n-\widehat{X}^-_n\Vert^r\right)^{1/r}\vee \EE\left(\Vert\, \widehat{m}_n-\widehat{X}_n\Vert^r\right)^{1/r}\right)\leq c(r)/\sqrt{N}
\end{equation}
\end{theo}
\proof
We have
$$
m_{n+1}-\widehat{X}^{-}_{n+1}=\displaystyle A ~(\widehat{m}_{n}-\widehat{X}_{n})+\frac{1}{\sqrt{N+1}}~R^{1/2}~\ZZ^{0}_{n+1}
$$
Using the decomposition
$$
\begin{array}{l}
\widehat{m}_{n}-\widehat{X}_{n}= \left(m_n-\widehat{X}^{-}_{n}\right)\\
\\
\displaystyle\hskip2cm+(K(p_n)-K(P_n))
~(Y_n-B\,  m_n) +K(P_n)(Y_n-B\,  m_n) -K(P_n)
~(Y_n-B  \widehat{X}^{-}_{n})\\
\\
\displaystyle\hskip3cm+\frac{1}{\sqrt{N+1}}~K(p_n)~R_0^{1/2}~\widehat{\ZZ}^{0}_n
\end{array}$$
we check that
\begin{equation}\label{stab-unif-w}
\begin{array}{l}
\widehat{m}_{n}-\widehat{X}_{n}= \left(\I-K(P_n)B\right)\left(m_n-\widehat{X}^{-}_{n}\right)\\
\\
\displaystyle\hskip1cm+\frac{1}{\sqrt{N}}~\sqrt{N}(K(p_n)-K(P_n))
~(Y_n-B\,  m_n) +\frac{1}{\sqrt{N+1}}~K(p_n)~R_0^{1/2}~\widehat{\ZZ}^{0}_n
\end{array}
\end{equation}
This implies that
\begin{equation}\label{stab-unif}
\left(m_{n+1}-\widehat{X}^{-}_{n+1}\right)= \Ea(P_n)\left(m_n-\widehat{X}^{-}_{n}\right)~+\frac{1}{\sqrt{N}}~\Omega_{n+1}
\quad \mbox{\rm with}\quad \Omega_{n+1}:=\Omega_{n}^{(1)}+\Omega^{(2)}_{n+1}
\end{equation}
In the above display $(\Omega_{n}^{(1)},\Omega^{(2)}_{n+1})$ stands for the remainder terms
\begin{eqnarray*}
 \Omega_{n}^{(1)}
&:=&A\sqrt{N}(K(p_n)-K(P_n))\left(B(X_n-m_n)+V_n\right)\\
\Omega^{(2)}_{n+1}&:=&\frac{1}{\sqrt{1+1/N}}~\left(AK(p_n)~R_0^{1/2}~\widehat{\ZZ}^{0}_n+R^{1/2}~\ZZ^{0}_{n+1}\right)
\end{eqnarray*}
Applying Theorem~\ref{th-intro-punif}, for any $r\geq 1$ and $N+1>d$ and we have the time-uniform estimates
$$
\sup_{n\geq 1}\EE\left(\Vert\, \Omega^{(2)}_{n}\Vert^r\right)\leq c(r)
$$
To take the final step, note that
$$
\begin{array}{l}
\displaystyle
m_{n+1}=A (\I-K(p_n)B)~ m_n+AK(p_n)
~Y_n +\frac{1}{\sqrt{N+1}}~\left(A~K(p_n)~R_0^{1/2}~\widehat{\ZZ}^{0}_n+R^{1/2}~\ZZ^{0}_{n+1}\right)
\end{array}$$
This yields
$$
\left(m_{n+1}-X_{n+1}\right)=\Ea(p_n)~ (m_n-X_n)+\Omega^{(3)}_{n+1}$$
with
$$
\Omega^{(3)}_{n+1}:=AK(p_n)V_n-W_{n+1} +\frac{1}{\sqrt{N+1}}~\left(A~K(p_n)~R_0^{1/2}~\widehat{\ZZ}^{0}_n+R^{1/2}~\ZZ^{0}_{n+1}\right)
$$
Applying Theorem~\ref{th-intro-punif}, for any $r\geq 1$ and $N+1>d$ and we have the time-uniform estimates
$$
\sup_{n\geq 1}\EE\left(\Vert\, \Omega^{(3)}_{n}\Vert^r\right)\leq c_1(r)
$$
Also note that
$$
\Ea(P)= S^{-1/2}~\overline{A}~S^{1/2}(\I+PS)^{-1}S^{-1/2}S^{1/2}=
S^{-1/2}~\left(\overline{A}~(\I+S^{1/2}PS^{1/2})^{-1}\right)S^{1/2}
$$
with
$$
\overline{A}:=S^{1/2}AS^{-1/2}\Longrightarrow
 \Vert \overline{A}~(\I+S^{1/2}PS^{1/2})^{-1}\Vert_2\leq  \Vert \overline{A}\Vert_2=\sqrt{\lambda_{1}(S^{1/2}AS^{-1}A^{\prime}S^{1/2})}
$$
This yields the almost sure exponential decays
$$
\Vert\Ea(p_n)\ldots \Ea(p_0)\Vert_2\leq  \left(\Vert S^{1/2}\Vert_2~\Vert S^{-1/2}\Vert_2\right) ~  \Vert  \overline{A}\Vert_2^{n+1}
$$
Thus using  corollary~\ref{cor-wpg}, for any $r\geq 1$ and $N+1>d$ we check the time-uniform estimates
$$
\sup_{n\geq 0}\EE\left(\Vert\, m_n-X_n\Vert^r\right)\leq c_2(r)\quad
\mbox{\rm and therefore}\quad
\sup_{n\geq 0}\EE\left(\Vert\, \Omega^{(1)}_{n}\Vert^r\right)\leq c_3(r)
$$
The end of the proof of (\ref{stab-unif-theo}) is now a direct consequence of
 the decomposition (\ref{stab-unif-w}) (\ref{stab-unif}) and the Floquet-type formula (\ref{form-Floquet}).
 This ends the proof of the theorem.\cqfd

\section*{Appendix}

\subsection*{Proof of (\ref{eq:Gn}),  (\ref{eq:gnB-ident1}), (\ref{square-root})
 and (\ref{def-Riccati-drift}).}\label{def-Riccati-drift-proof}

We recall the celebrated Sherman-Morrison-Woodbury matrix sum inversion identity
\begin{equation}\label{wood}
(M+UNV)^{-1}=M^{-1}-M^{-1}U(N^{-1}+VM^{-1}U)^{-1}VM^{-1},
\end{equation}
which is valid for any invertible matrices $(M,N)$ and any  conformable matrices $(U,V)$.

Choosing $M=\I$,  $U=P^{1/2}B^{\prime}$,
$V=BP^{1/2}$,  and $N=R_0^{-1}$ in (\ref{wood}) we check that
$$
\left(\I+P^{1/2}B^{\prime}R_0^{-1}BP^{1/2}\right)^{-1}=\I-P^{1/2}B^{\prime}(R_0+B  PB  ^{\prime} )^{-1}BP^{1/2}\Longleftrightarrow  (\ref{square-root})
$$

Applying the above to $M=\I$, $U=PB  ^{\prime}$, $N=-(B  PB  ^{\prime}+R_0 )^{-1}$ and $V=B  $ we check that
$$
\begin{array}{l}
(\I-K(P)B  )^{-1}=(\I-PB  ^{\prime}(B  PB  ^{\prime}+R_0 )^{-1}B  )^{-1}\\
\\
=\I-PB  ^{\prime}\left(-(B  PB  ^{\prime}+R_0 )+B  PB  ^{\prime}\right)^{-1}B  =
\I+PB  ^{\prime}(R_0)^{-1} B=\I+PS
\end{array}$$
This ends the proof of (\ref{eq:Gn}) and (\ref{def-Riccati-drift}).
Applying (\ref{wood}) to $M=P^{-1}$, $U=B^{\prime}$, $N=R_0^{-1}$ and $V=B$ we also check that
\begin{eqnarray*}
\widehat{A}(P)P&=&(\I-K(P)B)P=(\I+PS)^{-1}P=(P^{-1}+S)^{-1}\\
&=&P-PB^{\prime}\left(R_0+BPB^{\prime}\right)^{-1}BP
\end{eqnarray*}
Recalling that
$$
K(P)=PB^{\prime}\left(R_0+BPB^{\prime}\right)^{-1}
$$
we also check that
$$
\begin{array}{l}
\widehat{A}(P)P\widehat{A}(P)^{\prime}+\widehat{R}(P)\\
\\
=
(\I-K(P)B)P(\I-K(P)B)^{\prime}+K(P)\,R_0\,K(P)^{\prime}\\
\\
=\widehat{A}(P)P-\left(PB^{\prime}-K(P)\left(BPB^{\prime}+R_0\right)\right)K(P)^{\prime}=\widehat{A}(P)P
\end{array}$$
This ends the proof of (\ref{eq:gnB-ident1}). 
Note that
$$
(\I+PS)^{-1}P\leq  P\quad \mbox{\rm and}\quad(\I+PS)^{-1}P\leq  S^{-1}
$$
which implies that
$$
\Phi(P)\leq A ~S^{-1}~A ^{\prime}+R \quad \mbox{\rm and}\quad
R\leq \Phi(P)\leq A ~P~A ^{\prime}+R
$$
This ends the proof of (\ref{eq:Gn}) and (\ref{def-Riccati-drift}).
\cqfd
\subsection*{Proof of Lemma~\ref{lem-1-intro}}\label{lem-1-intro-proof}
We have the orthogonal decomposition
\begin{equation}\label{wish-dec}
 q(z+Z)=q(z)+R+\frac{1}{\sqrt{N}}~(\Ua+\Va(z))\end{equation}
with the centered and orthogonal $(d\times d)$-random symmetric matrices
$$
\Ua=\frac{1}{\sqrt{N}}
\sum_{1\leq i\leq N}\left(Z^i(Z^i)^{\prime}-R\right)\quad \mbox{and}\quad
\Va(z):=\frac{1}{\sqrt{N}}\sum_{1\leq i\leq N}\left(z^i(Z^i)^{\prime}+Z^i(z^i)^{\prime}\right)
$$

Notice that
$$
\begin{array}{l}
\displaystyle
zZ^{\prime}=\sum_{1\leq i\leq N} z^i(Z^i)^{\prime}=q^{1/2}(z)~\left(
\varsigma(z) ~\underline{Z}^{\prime}\right)~R^{1/2}
\quad \mbox{\rm and}\quad
ZZ^{\prime}=R^{1/2}(\underline{Z}\, \underline{Z}^{\prime}) R^{1/2}\\
\\
\displaystyle \mbox{\rm with}\quad
\varsigma(z):=q^{-1/2}(z)~z\quad \mbox{\rm and}\quad \underline{Z}:=R^{-1/2}~Z=[R^{-1/2}Z^1,\ldots,R^{-1/2}Z^N]\\
\\
\displaystyle\Longrightarrow
\frac{1}{N}~\varsigma(z)\varsigma(z)^{\prime}=q^{-1/2}(z)~\underbrace{\frac{1}{N}~zz^{\prime}}_{=q(z)}~q^{-1/2}(z)=q(\varsigma(z))=\I \\
\\
\displaystyle\Longleftrightarrow\forall 1\leq k,l\leq d\qquad
\frac{\varsigma(z)_k}{\sqrt{N}}~\left(\frac{\varsigma(z)_l}{\sqrt{N}}\right)^{\prime}=\left(\frac{\varsigma(z)^1_k}{\sqrt{N}},\ldots,\frac{\varsigma(z)_k^N}{\sqrt{N}}\right)\left(\begin{array}{c}
\frac{\varsigma(z)^1_l}{\sqrt{N}}
\\
\vdots\\
\frac{\varsigma(z)_l^N}{\sqrt{N}}
\end{array}\right)=1_{k=l}
\end{array}$$
We would like to show that the law of the random matrices
\begin{equation}\label{prop-equiv}
\left(\frac{1}{\sqrt{N}}~z ~\underline{Z}^{\prime}, \underline{Z}\,\underline{Z}^{\prime}\right)=\left(u ~\underline{Z}^{\prime}, \underline{Z}\,\underline{Z}^{\prime}\right)
\end{equation}
doesn't depends on the choice of the matrix $$
u:=\frac{1}{\sqrt{N}}~z \Longleftrightarrow
(u_1^{\prime},\ldots,u_d^{\prime}):=\left(\left(\frac{z_1}{\sqrt{N}}\right)^{\prime},\ldots,\left(\frac{z_d}{\sqrt{N}}\right)^{\prime}\right)$$
with $d$ orthogonal column vectors $u_k^{\prime}=\left(\begin{array}{c}
u_k^1
\\
\vdots\\
u_k^N
\end{array}\right)$ on the unit sphere of $\RR^N$.
That is, such that
$$
\forall 1\leq k,l\leq d\qquad u_ku_l^{\prime}=(u_k^1,\ldots,u_k^N)\left(\begin{array}{c}
u_l^1
\\
\vdots\\
u_l^N
\end{array}\right)=\sum_{1\leq i\leq N} u^i_k u^i_l=1_{k=l}
$$
Equivalently, in matrix forms
$$
\frac{1}{N}~zz^{\prime}=\left(\begin{array}{c}
\frac{z_1}{\sqrt{N}}
\\
\vdots\\
\frac{z_d}{\sqrt{N}}
\end{array}\right)\left(\left(\frac{z_1}{\sqrt{N}}\right)^{\prime},\ldots,\left(\frac{z_d}{\sqrt{N}}\right)^{\prime}\right)=\I=\left(\begin{array}{c}
u_1
\\
\vdots\\
u_d
\end{array}\right)\left(u_1^{\prime},\ldots,u_d^{\prime}\right)=uu^{\prime}
$$
Note that
$$
\pi_u=u^{\prime}u
$$
is the orthogonal projection from $\RR^N$ onto the vector space $\mbox{\rm Vect}(u^{\prime}_1,\ldots,u_d^{\prime})\subset \RR^N$ spanned by the orthonormal vectors $u^{\prime}_1,\ldots,u_d^{\prime}$.

Given another another  vector space $\mbox{\rm Vect}(v^{\prime}_1,\ldots,v_d^{\prime})\subset \RR^N$ spanned by $d$ orthonormal vectors $v^{\prime}_1,\ldots,v_d^{\prime}$, there exists a rotation matrix $R^{\prime}$ such that $v^{\prime}_k=R^{\prime}u^{\prime}_k$. By construction, we have
$$
v \underline{Z}^{\prime}=\left(\begin{array}{c}
v_1
\\
\vdots\\
v_d
\end{array}\right)\left( \underline{Z}_1^{\prime},\ldots, \underline{Z}_d^{\prime}\right)=\left(\begin{array}{c}
u_1
\\
\vdots\\
u_d
\end{array}\right)\left( R\underline{Z}_1^{\prime},\ldots, R\underline{Z}_d^{\prime}\right)
$$
In addition, we have
$$
\underline{Z}R^{\prime}R\underline{Z}^{\prime}=
\underline{Z}\,\underline{Z}^{\prime}
\quad \mbox{\rm and}\quad
R\underline{Z}^{\prime}=\left(R\underline{Z}_1^{\prime},\ldots, R\underline{Z}_d^{\prime}\right)\stackrel{law}{=}\underline{Z}^{\prime}=
\left(\underline{Z}_1^{\prime},\ldots, \underline{Z}_d^{\prime}\right)
$$
By the rotational invariance of Gaussian random vectors, we conclude that
$$
\left(v ~\underline{Z}^{\prime}, \underline{Z}\,\underline{Z}^{\prime}\right)=\left(u ~R\underline{Z}^{\prime}, \underline{Z}\,R^{\prime}R\,\underline{Z}^{\prime}\right)\stackrel{law}{=}\left(u ~\underline{Z}^{\prime}, \underline{Z}\,\,\underline{Z}^{\prime}\right)
$$
Choosing
$$
\begin{array}{l}
\displaystyle
\forall 1\leq k\leq d\qquad\forall 1\leq i\leq N\qquad
u^i_k=1_{i=k}\\
\\
\Longrightarrow \forall 1\leq k,l\leq d\qquad
\left(u ~\underline{Z}^{\prime}\right)_{k}^l=\sum_{1\leq i\leq N} u^i_k~\underline{Z}^i_l=
\underline{Z}^k_l=(\underline{Z}^{\prime})_k^l
\end{array}$$
This ends the proof of (\ref{wish-dec-intro}).
\cqfd

\subsection*{Proof of Proposition~\ref{var-ncW-prop}.}\label{var-ncW-proof}

After some elementary computations, for any $A\in\Ma_d$ we check that
 $$
 \EE(\HH^NA~\HH^N)=A^{\prime}+\tr(A)~\I \quad \mbox{\rm and}\quad
 \EE(\HH^N A\GG)=0
 $$
 as well as
 $$
 \EE(\GG A\GG)=A^{\prime}\quad \mbox{\rm and}\quad
  \EE(\GG A\GG^{\prime})=\tr(A)~\I
 $$
 This implies that
  $$
  \begin{array}{l}
 \EE\left(\left(q^{1/2}(z)~ \GG^{\prime}~R^{1/2}+R^{1/2}~ \GG~q^{1/2}(z)\right)^2\right)\\
 \\
 = q^{1/2}(z)~\EE\left( \GG~R^{1/2}q^{1/2}(z)~\GG\right)R^{1/2}+R^{1/2} \EE\left(\GG~q^{1/2}(z)R^{1/2} \GG\right)~~q^{1/2}(z)\\
 \\
\hskip3cm +q^{1/2}(z)~\EE\left( \GG^{\prime}~R~ \GG~\right)q^{1/2}(z)~+\EE\left(R^{1/2} ~\GG~q(z)~ \GG^{\prime}\right)R^{1/2}\\
\\
=q(z)R+Rq(z)+\tr(R) q(z)+\tr(q(z))R
 \end{array}
 $$
 from which we check that
 $$
 \EE\left(\Delta^N_q(z)^2\right)=R^{2}+\tr(R)~R+
q(z)R+Rq(z)+\tr(R) q(z)+\tr(q(z))R
 $$
 This implies that
 $$
   \begin{array}{l}
\EE\left( \left( q(z+Z)-(q(z)+R)\right)^2\right)=\EE\left(q(z+Z)^2\right)-(q(z)+R)^2\\
\\
\displaystyle=\frac{1}{N}~\left(R^{2}+\tr(R)~R+
q(z)R+Rq(z)+\tr(R) q(z)+\tr(q(z))R\right)
  \end{array}
 $$
 This ends the proof of (\ref{var-ncW}). The proof of the proposition is now completed.\cqfd

\subsection*{Proof of Lemma~\ref{lem-2-intro}}\label{lem-2-intro-proof}
We fix some parameter $N\geq 1$ and we let $\un$ be the $(N+1)$ column vector with unit entries, $\I$ the $(N+1)\times(N+1)$ identity matrix and $\JJ=\un\un^{\prime}$ the $(N+1)\times(N+1)$ matrix with unit entries, where $(\cdot)'$ denotes the transpose operator. We also let $\epsilon$ be  the $(N+1)\times(N+1)$ matrix given by
$$
\epsilon= \I-{\JJ}/{(N+1)}
\quad \Longrightarrow \quad
\epsilon^2=\epsilon=\epsilon^{\prime}\quad \mbox{\rm and}\quad\sum_{1\leq j\leq N+1}\epsilon^i_j=0
$$
Further, denote by $\OO$ the $((N+1)\times(N+1))$-Helmert matrix whose first row $\OO_1$ is defined  by
$$
 \OO_1=\un^{\prime}/\sqrt{N+1}.
$$
and whose $i$-th row for $2\leq i\leq (N+1)$ is given by
$$
\begin{array}{l}
\OO_i^j:=\overline{\OO}_{i-1}^j:=\left\{
\begin{array}{ccl}
1/\sqrt{i(i-1)}&\mbox{\rm if}& 1\leq j<i\\
&&\\
-(i-1)/\sqrt{i(i-1)}&\mbox{\rm if}& j=i\\
&&\\
0&\mbox{\rm if}& i<j\leq (N+1).
\end{array}\right.
\end{array}
$$
This yields the matrix decomposition
$$
\OO=\left(\begin{array}{c}
\OO_1\\
\overline{\OO}\end{array}\right)\quad \mbox{\rm with}\quad \overline{\OO}=\left[\begin{array}{c}
\overline{\OO}_1\\
\vdots\\
\overline{\OO}_N
\end{array}\right]=\left[\overline{\OO}^1,\ldots, \overline{\OO}^{N+1}\right]\in \RR^{N\times (N+1)}.
$$
We check that $\OO$ is an orthogonal matrix and we have the square root formula
\begin{equation}\label{H-epsi}
\begin{array}{l}
\displaystyle\OO\,\OO^{\prime}=\I=\OO^{\prime}\,\OO=\frac{1}{N+1}~\JJ+\overline{\OO}^{\prime}~\overline{\OO} \Longrightarrow \overline{\OO}^{\prime}~\overline{\OO}=\epsilon.
\end{array}
\end{equation}
We set
$$(Z^0,Z):=(Z^0,(Z^1,\ldots,Z^N))=\Za~\OO^{\prime} \in\Ma_{d\times (N+1)}\quad \mbox{\rm and}\quad
 \overline{x}:=x~\overline{\OO}^{\prime}\in \RR^{d\times N}
$$
By construction,  $Z^i$ are $(N+1)$ iid centered Gaussian  $d$-column vectors
 with covariance matrix $R$ and we have
\begin{eqnarray*}
p(x)
&=&\frac{1}{N}~x~ \overline{\OO}^{\prime}~\overline{\OO}~x^{\prime}=\frac{1}{N}~ \overline{x}~ \overline{x}^{\prime}=q( \overline{x})
\end{eqnarray*}
The last assertion comes from the matrix formula
\begin{equation}
x^{\epsilon}:=x~\epsilon=(x-M(x))\label{ref-1-U}
\end{equation}
Observe that
\begin{eqnarray}
(Z^0,Z)=(\Za~\OO^{\prime}_1,\Za~\overline{\OO}^{\prime})&\Longleftrightarrow& Z^0=\sqrt{N+1}~m(\Za)\quad \mbox{\rm and}\quad Z=\Za~\overline{\OO}^{\prime}\label{ref-1Weps}
\end{eqnarray}
Using  \eqref{ref-1Weps} we check that
$$
\Za^{\epsilon}(\Za^{\epsilon})^{\prime}=\Za~\overline{\OO}^{\prime}~\overline{\OO}~\Za^{\prime}=Z~Z^{\prime}
\quad\mbox{\rm and}\quad
m(x+\Za)=m(x)+\frac{1}{\sqrt{N+1}}~Z^0
$$
Applying \eqref{ref-1-U}, we have
$$
p(x+\Za)= \frac{1}{N}~(x^{\epsilon}+\Za^{\epsilon})(x^{\epsilon}+\Za^{\epsilon})^{\prime}
$$
This implies that
$$
p(x+\Za)=p(x)+\frac{1}{N}~Z~Z^{\prime}+\frac{1}{N}~\left(\overline{x}~Z^{\prime}+Z~\overline{x}^{\prime}\right)=q(\overline{x}+Z)
$$
The end of the proof of Lemma~\ref{lem-2-intro} is now a direct consequence of Lemma~\ref{lem-1-intro}.
\cqfd

\subsection*{Proof of corollary~\ref{cor-wpg}}\label{cor-wpg-proof}
Recalling that
$
(\I-K(P)B)=(\I+PS)^{-1}
$,
we check that
$$
\begin{array}{l}
\displaystyle
(K(p_n)-K(P_n))B=(\I+P_nS)^{-1}-(\I+p_nS)^{-1}\\
\\
\displaystyle=
~(\I+P_nS)^{-1}\left(p_n-P_n\right)S(\I+p_nS)^{-1}=
~S^{-1}(S^{-1}+P_n)^{-1}\left(p_n-P_n\right)(S^{-1}+p_n)^{-1}
\end{array}
$$
This implies that
$$
\begin{array}{l}
\displaystyle
\Vert (K(p_n)-K(P_n))B\Vert_{\tiny F}\\
\\
\displaystyle\leq  \Vert S^{-1}\Vert_{\tiny F}~~\Vert (S^{-1}+P_n)^{-1}\Vert_{\tiny F}~~
\Vert (S^{-1}+p_n)^{-1}\Vert_{\tiny F}~~\Vert p_n-P_n\Vert_{\tiny F}\\
\\
\displaystyle\leq \lambda_{1}(S)~\tr(S)~ \Vert S^{-1}\Vert_{\tiny F}~~\Vert p_n-P_n\Vert_{\tiny F}
\end{array}
$$
The last assertion comes from the fact that
\begin{eqnarray*}
(S^{-1}+P)^{-2}&=&
(S^{-1}+P)^{-1/2}(S^{-1}+P)^{-1}
~(S^{-1}+P)^{-1/2}\\
&\leq&
(S^{-1}+P)^{-1/2}S
~(S^{-1}+P)^{-1/2}
\leq \lambda_{1}(S)~(S^{-1}+P)^{-1}\leq  \lambda_{1}(S)~S
\end{eqnarray*}
We conclude that
$$
\Vert K(p_n)-K(P_n)\Vert_{\tiny F}\leq\lambda_{d}(BB^{\prime})^{-1/2}
 \lambda_{1}(S)~\tr(S)~ \Vert S^{-1}\Vert_{\tiny F}~~\Vert p_n-P_n\Vert_{\tiny F}
$$
Finally observe that
$$
\widehat{p}_n-\widehat{P}_n=\left(\widehat{A}(p_n)~ p_n-\widehat{A}(P_n)~ P_n\right)+\frac{1}{\sqrt{N}}~ \widehat{\Delta}^N_{n}
$$
with the mapping $\widehat{A}$ defined in (\ref{eq:Gn}).
On the other hand, we have
\begin{eqnarray*}
\widehat{A}(p_n)~ p_n-\widehat{A}(P_n)~ P_n&=&(K(P_n)-K(p_n))B  P_n+(\I-K(p_n)B  )(p_n-P_n)\\
&=&(K(P_n)-K(p_n))B  P_n+(\I+p_nS)^{-1}(p_n-P_n)
\end{eqnarray*}
This implies that
$$
\begin{array}{l}
\Vert \widehat{A}(p_n)~ p_n-\widehat{A}(P_n)~ P_n\Vert_{\tiny F}\\
\\
\leq \Vert (K(P_n)-K(p_n))B\Vert_{\tiny F}~\Vert  P_n \Vert_{\tiny F}+
\Vert S^{-1}\Vert_{\tiny F}
\Vert(S^{-1}+p_n)^{-1}\Vert_{\tiny F}~\Vert p_n-P_n\Vert_{\tiny F}
\end{array}
$$
Recalling that
$$
 P_n\leq  AS^{-1}A^{\prime}+R\quad \mbox{\rm and}\quad
\Vert(S^{-1}+p_n)^{-1}\Vert_{\tiny F}\leq (\lambda_{1}(S)~\tr(S))^{1/2}
$$
there exists some finite constant $c>0$ such that for any $n\geq 1$ we have
$$
\Vert \widehat{A}(p_n)~ p_n-\widehat{A}(P_n)~ P_n\Vert_{\tiny F}\leq c~\Vert p_n-P_n\Vert_{\tiny F}
$$
The end of the proof of (\ref{unif-cov-2}) is now a direct consequence of the uniform estimates (\ref{unif-cov}). This ends the proof of Corollary~\ref{cor-wpg}.\cqfd

\section*{No conflict of interest}

On behalf of all authors, the corresponding author states that there is no conflict of interest.

\end{document}